\documentclass{article}

\usepackage{amssymb}
\usepackage{amsmath}
\usepackage{mathrsfs}
\usepackage{amsthm}

\theoremstyle{plain}
\newtheorem{thm}{Theorem}[section]

\newtheorem{prp}[thm]{Proposition}
\newtheorem{cor}[thm]{Corollary}

\theoremstyle{definition}
\newtheorem{dfn}[thm]{Definition}

\newtheorem*{Ter}{Terminology}
\newtheorem*{Not}{Notation}

\theoremstyle{remark}
\newtheorem{rem}[thm]{Remark}
\newtheorem{exa}[thm]{Example}
\newtheorem{Q}[thm]{Question}

\def\({{\rm (}}
\def\){\/{\rm )}} 
\def\le{\leqslant}
\def\ge{\geqslant}
\def\]{\mathopen]}
\def\[{\mathclose[}

\mathchardef\hook="312C

\def\projtens{\mathop{\otimes}\limits^\wedge}
\def\injtens{\mathop{\otimes}\limits^\vee}
\def\concat{\mathord{{}^\smallfrown}}
\def\imp{\Rightarrow}
\def\ssi{\Leftrightarrow}
\def\card#1{\# #1}

\def\trace{\mathop{\rm\rm tr}\nolimits} 
\def\id{{\rm Id}}
\def\Alpha{{\rm A}}
\def\P#1{{\rm P}_{\mkern-6mu#1}}
\def\Pw#1{\mathscr{P}_{\mkern-6mu#1}}
\def\e{\mkern1mu{\rm e}\mkern1mu}
\def\iu{\mkern1mu{\rm i}\mkern1mu}
\def\Ref#1{\hbox{$(\ref{#1})$}}

\def\Z{\mathbb{Z}}
\def\U{\mathbb{S}}
\def\T{\mathbb{T}}
\def\N{\mathbb{N}}
\def\R{\mathbb{R}}
\def\C{\mathbb{C}}

\def\row{R}
\def\col{C}
\title{Cycles and $1$-unconditional matrices}
\author{Stefan Neuwirth}
\begin{document}
\parindent0pt
\maketitle

\begin{abstract}
\noindent
We characterize the $1$-unconditional subsets
$(\e_{rc})_{(r,c)\in I}$ of the set of elementary matrices
in the Schatten-von-Neumann class $S^p$. The set of couples $I$ must
be the set of edges of a bipartite graph without cycles of even length
$4\leqslant l\leqslant p$ if $p$ is an even integer, and without
cycles at all if $p$ is a positive real number that is not an even
integer. In the latter case, $I$ is even a Varopoulos set of
$V$-interpolation of constant $1$. We also study the metric unconditional
approximation property for the space $S^p_I$ spanned by
$(\e_{rc})_{(r,c)\in I}$ in $S^p$.
\end{abstract}

\section{Introduction}

The starting point for this investigation has been the following
isometric question on the Schatten-von-Neumann class $S^p$.
\begin{Q}\label{q1}
  Which matrix coefficients of an operator $x\in S^p$ must vanish so
  that the norm of $x$ does not depend on the argument, or on the
  sign, of the remaining nonzero matrix coefficients~?
\end{Q}
Let $\col$ be the set of columns and $\row$ be the set of rows for
coordinates in the matrix. Let $I\subseteq\row\times\col$ be the set of
matrix coordinates of the nonzero matrix coefficients of $x$.
Question~\ref{q1} describes the notion of a complex, or real,
$1$-unconditional basic sequence $(\e_{rc})_{(r,c)\in I}$ of
elementary matrices in $S^p$. \medskip

By a convexity argument, Question~\ref{q1} is
equivalent to the following question on Schur multiplication.
\begin{Q}\label{q2}
  Which matrix coefficients of an operator $x\in S^p$ must vanish so
  that for all matrices $\varphi$ of complex, or real, numbers
  $\|\varphi*x\|\leqslant\sup|\varphi_{rc}|~\|x\|$, where $*$ is
  the entrywise (or Schur or Hadamard) product given by
  $(\varphi*x)_{rc}=\varphi_{rc}x_{rc}$~?
\end{Q}

In the case $p=\infty$, Grothendieck's inequality yields an estimation
for the norm of Schur multiplication by $\varphi$ in terms of the
projective tensor product $\ell_\col^\infty\projtens\ell_\row^\infty$:
it is equivalent to the supremum of the norm of the tensors whose
coefficient matrices are finite submatrices of $\varphi$. In the
framework of tensor algebras over discrete spaces, Question~\ref{q2}
turns out to describe as well the isometric counterpart to Varopoulos'
$V$-Sidon sets as well as to his sets of $V$-interpolation. The
following isometric question has however a different answer.
\begin{Q}\label{q3}
  Which coefficients of a tensor
  $u\in\ell_\col^\infty\projtens\ell_\row^\infty$ must vanish so that
  the norm of $u$ is the maximal modulus of its coefficients~? 
\end{Q}

In our answer to Question~\ref{q2}, $S^p$ and Schur multiplication
are treated as a non commutative analogue to $L^p$ and convolution.
The main step is a careful study of the Schatten-von-Neumann norm
$\|x\|=\bigl(\trace(x^*x)^{p/2}\bigr)^{1/p}$ for $p$ an even integer.
The rule of matrix multiplication provides an expression for this norm
as a series in the matrix coefficients of $x$ and their complex
conjugate, indexed by the $p$-uples $(v_1,v_2,\dots,v_p)$ satisfying
$(v_{2i-1},v_{2i}),(v_{2i+1},v_{2i})\in I$. Those are best understood
as closed walks of length $p$ on the bipartite graph $G$ canonically
associated to $I$: its vertex classes are $\col$ and $\row$ and its
edges are given by the couples in $I$. A structure theorem for closed
walks and a detailed study of the particular case in which $G$ is a cycle yield
the two following theorems.
\begin{thm}\label{th4}
  Let $p\in\]0,\infty]\setminus\{2,4,6,\dots\}$. If the
  sequence of elementary matrices $(\e_{rc})_{(r,c)\in I}$ is a real
  $1$-unconditional basic sequence in $S^p$, then $G$ contains no
  cycle. In this case, $I$ is even a set of $V$-interpolation with
  constant $1$: every sequence $\varphi\in\ell^\infty_I$ may be
  interpolated by a tensor
  $u\in\ell_\col^\infty\projtens\ell_\row^\infty$ such that
  $\|u\|=\|\varphi\|$.
\end{thm}
\begin{thm}\label{th5}
  Let $p\in\{2,4,6,\dots\}$. The sequence
  $(\e_{rc})_{(r,c)\in I}$ is a complex, or real, $1$-unconditional
  basic sequence in $S^p$ if and only if $G$ contains no cycle of
  length $4,6,\dots,p$.
\end{thm}
If $C$ and $R$ are finite, extremal graphs without cycles of given
lengths remain an ongoing area of research in Graph theory.\medskip

These results provide a complete description of the situation in which
$(\e_{rc})_{(r,c)\in I}$ is a $1$-un\-con\-di\-tio\-nal basis of the space
$S^p_I$ it spans in $S^p$. If this is not the case, $S^p_I$ might
still admit some other $1$-unconditional basis. This leads to the
following more general question. 
\begin{Q}\label{q4}
  For which sets $I$ does $S^p_I$ admit some kind of almost
  $1$-un\-con\-di\-tional Schauder decomposition~?
\end{Q}
The metric unconditional approximation property $(umap)$ provides a formal
definition for the object of Question~\ref{q4}. We obtain the following results.
\begin{thm}
  Let $p\in[1,\infty]\setminus\{2,4,6,\dots\}$. If $S^p_I$ has real
  $(umap)$, then the distance of any two vertices $c\in\col$ and
  $r\in\row$ is asymptotically infinite in $G$: their distance becomes
  arbitrarily large by deleting a finite number of edges from $G$.
\end{thm}
\begin{thm}
  Let $p\in\{2,4,6,\dots\}$. The space $S^p_I$ has complex, or real,
  $(umap)$ if and only if any two vertices at distance $2j+1\le p/2$
  are asymptotically at distance at least $p-2j+1$.
\end{thm}
 
We now turn to a detailed description of this article. In
Section~\ref{ss2}, we 
provide tools for the computation of Schur
multiplier norms. Section~\ref{ss3} characterizes idempotent Schur
multipliers and $0,1$-tensors in $\ell_\col^\infty\projtens\ell_\row^\infty$
of norm one. In Section~\ref{ss4}, we define the complex and real
unconditional constants of basic sequences of elementary matrices and
show that they are not equal in general. Section~\ref{V} looks back on
Varopoulos' results about tensor algebras over discrete spaces.
Section~\ref{ss6} puts the connection between $p$-trace norm and
closed walks of length $p$ in the concrete form of closed walk
relations. In Section~\ref{ss7}, we compute the norm of relative Schur
multipliers by signs in the case that $G$ is a cycle, and estimate the
corresponding unconditional constants. Section~\ref{ss8} is dedicated
to a proof of Th.~\ref{th4} and an answer to Question~\ref{q3}.
Section~\ref{ss9} establishes Th.~\ref{th5}. In Section~\ref{ss10},
we study the metric unconditional approximation property for spaces
$S^p_I$. The final section provides four kinds of examples: sets
obtained by a transfer of $n$-independent subsets of a discrete abelian
group, Hankel sets, Steiner systems and Tits' generalized polygons.

\begin{Ter}
  $\col$ is the set of \textit{columns} and $\row$ is the set of
  \textit{rows}, both finite or countable and if necessary indexed by
  natural numbers.  $V$, the set of \textit{vertices}, is their
  disjoint union $\col\amalg\row$. An \textit{edge} on $V$ is a pair
  $\{v,w\}\subseteq V$. A \textit{graph} on $V$ is given by a set of
  edges $E$. A \textit{bipartite} graph on $V$ with vertex classes
  $\col$ and $\row$ has only edges $\{r,c\}$ such that $c\in\col$ and
  $r\in\row$ and may therefore be given alternatively by the set of
  couples $I=\bigl\{(r,c)\in\row\times\col:\{r,c\}\in E\bigr\}$. Two
  graphs are \textit{disjoint} if so are the sets of vertices of their
  edges. $I$ is a \textit{row section} if $(r,c),(r,c')\in I\imp
  c=c'$, and a \textit{column section} if $(r,c),(r',c)\in I\imp
  r=r'$.\medskip

  A \textit{walk} of length $s\geqslant0$ in a graph is a sequence
  $(v_0,\dots,v_s)$ of $s+1$ vertices such that
  $\{v_0,v_1\},\dots,\{v_{s-1},v_s\}$ are edges of the graph. A walk
  is a \textit{path} if its vertices are pairwise distinct. The
  \textit{distance} of two vertices in a graph is the minimal length
  of a path in the graph that joins the two vertices; it is infinite
  if no such path exists. A \textit{closed} walk of length $p$ in a
  graph is a sequence $(v_1,\dots,v_p)$ of $p$ vertices such that
  $\{v_1,v_2\},\dots,\{v_{p-1},v_p\},\{v_p,v_1\}$ are edges of the
  graph. Note that $p$ is necessarily even if the graph is bipartite.
  A closed walk is a \textit{cycle} if its vertices are pairwise
  distinct. We take the convention that the first vertex of a closed
  walk or a cycle on $V=\col\amalg\row$ is a column vertex:
  $v_1\in\col$. We shall identify a path and a cycle with
  their set of edges $\{r,c\}$ or the corresponding set of couples
  $(r,c)$.\medskip

  A bipartite graph on $V$ is a \textit{tree} if its vertices may be
  indexed by a set $W$ of finite words over some set of letters $A$ in
  the following way:
\begin{itemize}
\item $\emptyset\in W$ and every beginning of a word in $W$ is also in
  $W$: if $w\in W\setminus\{\emptyset\}$, then $w$ is the
  concatenation $w'\concat a$ of a word $w'\in W$ with a letter $a\in
  A$;
\item words of even ({\it vs.}~odd) length index row ({\it vs.}~column)
  vertices;
\item a pair of vertices is
an edge exactly if the length of their indices differs by one.
\end{itemize}
A \textit{forest} is a union of pairwise disjoint trees; equivalently, it is a cycle free
graph.
\end{Ter}

\begin{Not}
$\T=\{z\in\C:|z|=1\}$. Given an index set
$I$ and $q\in I$, $\e_q$ is the sequence defined on $I$ as the
indicator function $\chi_{\{q\}}$ of the singleton $\{q\}$.\medskip

Let $I=\row\times\col$ and $q=(r,c)$. Then $\e_q=\e_{rc}$ 
is the elementary matrix identified with the operator from
$\ell^2_\col$ to $\ell^2_\row$ that maps $\e_c$ on $\e_r$ and all
other basis vectors on $0$. The \textit{matrix coefficient} at
coordinate $q$ of an operator $x$ from $\ell^2_\col$ to $\ell^2_\row$
is $x_q=\trace\e_q^*x$ and its \textit{matrix representation} is
$(x_q)_{q\in\row\times\col}=\sum_{q\in\row\times
  \col}x_q\e_q$.\medskip

The \textit{Schatten-von-Neumann class} $S^p$, $0<p<\infty$, is the space of those
compact operators $x$ from $\ell^2_\col$ to $\ell^2_\row$ such that
$\|x\|_p^p=\trace|x|^p=\trace(x^*x)^{p/2}<\infty$. $S^\infty$ is the space of
compact operators with the operator norm. $S^p$ is a quasi-normed space, and a
Banach space if $p\geqslant1$. 
The unit ball of a Banach space $X$ is denoted by $B_X$.\medskip

For $I\subseteq\row\times\col$, the \textit{entry space} $S^p_I$ is the space
of those $x\in S^p$ whose \textit{support}
$\{q\in\row\times\col:x_q\ne0\}$ is a subset of $I$. $S^p_I$ is also the closed subspace of $S^p$
spanned by $(\e_q)_{q\in I}$.\medskip

The $S^p$-valued Schatten-von-Neumann class $S^p(S^p)$ is the space of those
compact operators $x$ from $\ell^2_\col$ to $\ell^2_\row(S^p)$ such
that $\|x\|_p^p=\trace(\trace|x|^p)<\infty$, where the inner trace is
the $S^p$-valued analogue of the usual trace $x$: such operators have
an $S^p$-valued matrix representation.  $S^p(S^p)$ can be identified
with the space of compact operators $x$ from $\ell^2_\col(\ell_2)=\ell_2\otimes_2\ell^2_\col$ to
$\ell^2_\row(\ell_2)=\ell_2\otimes_2\ell^2_\row$ such that
$\|x\|_p^p=\trace\otimes\trace|x|^p<\infty$; the matrix coefficient of
$x$ at $q$ is then
$x_q=(\id_{S^p}\otimes\trace)\bigl((\id_{\ell_2}\otimes\e_q^*)x\bigr)$ and its
matrix representation is $\sum_{q\in\row\times \col}x_q\otimes\e_q$.\medskip

A \textit{relative Schur multiplier} on $S^p_I$ is a sequence
$\varphi=(\varphi_q)_{q\in I}\in\C^I$ such that the operator
$M_\varphi$ defined by $\e_q\mapsto\varphi_q\e_q$ for $q\in I$ is
bounded on $S^p_I$. The Schur multiplier $\varphi$ is furthermore
\textit{completely bounded} (c.b.~for short) on $S^p_I$ if
$\id_{S^p}\otimes M_\varphi$, that is the operator defined by
$x_q\e_q\mapsto\varphi_qx_q\e_q$ for $x_q\in S^p$ and $q\in I$, is
bounded on $S^p_I(S^p)$ (see \cite[Lemma~1.7]{pi98}). Note that
$\varphi$ is a Schur multiplier on $S^\infty$ if and only if it is a
Schur multiplier on the space of bounded operators from $\ell^2_\col$
to $\ell^2_\row$ and that it automatically is c.b.\ 
\cite[Th.~5.1]{pi01}.\medskip

Let $X,Y$ be Banach spaces and $u\in X\otimes Y$. Its \textit{projective}
tensor norm is given by 
$$
\|u\|_{X\projtens Y}=\inf\biggl\{\sum_{j=1}^n\|x_j\|\,\|y_j\|:u=\sum_{j=1}^nx_j\otimes
y_j\biggr\}
$$
and $X\projtens Y$ is the completion of $X\otimes Y$ with
respect to this norm. Note that 
$\ell_\infty^n\projtens \ell_\infty^m\subset
c_0\projtens c_0\subset\ell_\infty\projtens \ell_\infty$. 
The \textit{injective} tensor norm
of $u$ is given by
$$
\|u\|_{X\injtens Y} =
\sup_{(\xi,\eta)\in B_{X^*}\times B_{Y^*}}|\langle\xi\otimes\eta,u\rangle|
$$
and $X\injtens Y$ is the completion of $X\otimes Y$ with
respect to this norm. Note that if $X$ and $Y$ are both finite dimensional, then
$$
(X\injtens Y)^*=X^*\projtens Y^*\quad{\rm
  and}\quad(X\projtens Y)^*=X^*\injtens Y^*~;
$$
further
$(c_0\projtens c_0)^*=\ell_1\injtens\ell_1$
and $\ell_\infty\projtens \ell_\infty$ is a closed $w^*$-dense
subspace of $(\ell_1\injtens\ell_1)^*$. 
If $X$ is a sequence space on $\col$ and $Y$ is a sequence space on $\row$,
then the \textit{coefficient} of the tensor $u$ at $(r,c)$ is
$\langle\e_c\otimes\e_r,u\rangle$. Its \textit{support} is the set of
coordinates $(r,c)$ of its nonvanishing coefficients.
\end{Not}

\section{Relative Schur multipliers}\label{ss2}

The following proposition is a straightforward consequence of
\cite{pps89}.

\begin{prp}\label{pinf}

  Let $I\subseteq\row\times\col$ and $\varphi$ be a Schur multiplier
  on $S^\infty_I$ with norm $D$. Then $\varphi$ is also a c.b.~Schur multiplier
  on $S^p_I$ for every $p\in\]0,\infty]$, with c.b.~norm bounded by $D$.

\end{prp}

\begin{proof}
  
  We may assume $D=1$. Let $\row'\times\col'$ be any finite subset of
  $\row\times\col$. By \cite[Th.~3.2]{pps89}, there exist vectors
  $w_c$ and $v_r$ of norm at most $1$ in a Hilbert space
  $H$ such that $\varphi_{rc}=\langle w_c,v_r\rangle$ for every
  $(r,c)\in I\cap\row'\times\col'$. If we define $W\colon\ell^2_{C'}\to\ell^2_{C'}(H)$ and
  $V\colon\ell^2_{R'}\to\ell^2_{R'}(H)$ by
  $W\zeta=(\zeta_cw_c)_{c\in C'}$ and $V\eta=(\eta_rv_r)_{r\in R'}$,
  then $V$ and $W$ have norm at most $1$, and the proposition follows from the factorization
\begin{equation*}
M_\varphi x=V^*(x\otimes\id_H)W
\end{equation*}
for every $x$ with support in $I\cap\row'\times\col'$.
\end{proof}

\begin{rem}\label{rem} 
  \'Eric Ricard showed us an elementary proof
  that a Schur multiplier on $S^\infty_I$ is automatically
  c.b., included here by his kind permission. A Schur multiplier
  $\varphi$ is bounded on $S^\infty_I$ by a constant $D$ if and only if
\begin{equation}\label{b}
\forall\,\xi\in
B_{S^\infty_I}~\forall\,\eta\in B_{\ell^2_\row}~\forall\,\zeta\in
B_{\ell^2_\col}\quad\biggl|\sum_{(r,c)\in
  I}\eta_r\varphi_{rc}\xi_{rc}\zeta_c\biggr|\leqslant D~.
\end{equation}
It is furthermore completely bounded on $S^\infty_I$ by $D$ if
\begin{equation}\label{cb}
\forall\,x\in
B_{S^\infty_I(S^\infty)}~\forall\,y\in B_{\ell^2_\row(\ell_2)}~\forall\,z\in
B_{\ell^2_\col(\ell_2)}\quad\biggl|\sum_{(r,c)\in I}\varphi_{rc}\langle
y_r,x_{rc}z_c\rangle\biggr|\leqslant D~.
\end{equation}
Now, given $x,y,z$ as quantified in Ineq.~\Ref{cb}, let $\xi_{rc}=\langle
y_r/\|y_r\|,x_{rc}z_c/\|z_c\|\rangle$, $\eta_r=\|y_r\|_{\ell_2}$,
$\zeta_c=\|z_c\|_{\ell_2}$. Then
$\|\eta\|_{\ell^2_\row},\|\zeta\|_{\ell^2_\col}\leqslant1$ and
\begin{eqnarray*}
  \|\xi\|&=&\sup\biggl\{\biggl| \sum_{(r,c)\in I} \langle \alpha_ry_r/\|y_r\|_{\ell_2},
  x_{rc}\beta_c z_c/\|z_c\|_{\ell_2} \rangle \biggr| : \alpha\in B_{\ell^2_\row},\,
  \beta\in B_{\ell^2_\col} \biggr\}\\
&\leqslant& \|x\| 
\sup_{\alpha\in B_{\ell^2_\row}}
\bigl\|\bigl(\alpha_ry_r/\|y_r\|_{\ell_2}\bigr)\bigr\|_{\ell^2_\row(\ell_2)}
\sup_{\beta\in B_{\ell^2_\col}}
\bigl\|\bigl(\beta_cz_c/\|z_c\|_{\ell_2}\bigr)\bigr\|_{\ell^2_\col(\ell_2)}\leqslant 1~, 
\end{eqnarray*}
so that Ineq.~\Ref{b} implies Ineq.~\Ref{cb}.
\end{rem}

The fact that the canonical basis of an $\ell^2$ space is
$1$-unconditional yields that Schatten-von-Neumann norms are
\textit{matrix unconditional} in the terminology of \cite{sc78}:
\begin{equation}
  \label{mu}
  \forall\,\zeta\in\T^\col~\forall\,\eta\in\T^\row\quad\Bigl\|\sum_{(r,c)\in\row\times\col}\zeta_c\eta_ra_{rc}\e_{rc}\Bigr\|_p=\Bigl\|\sum_{(r,c)\in\row\times\col}
  a_{rc}\e_{rc}\Bigr\|_p
\end{equation}
for every finitely supported sequence of complex or $S^p$-valued
coefficients $a_{rc}$. This shows that if $\zeta\in\ell^\infty_\col$,
$\eta\in\ell^\infty_\row$, then the c.b.~norm of the \textit{elementary} Schur
multiplier $(\zeta_c\eta_r)_{(r,c)\in I}$ is bounded by
$\|\zeta\|_{\ell^\infty_\col}\|\eta\|_{\ell^\infty_\row}$ on every
$S^p_I$.\medskip

The following proposition relates Fourier multipliers to Herz-Schur
multipliers in the fashion of \cite[Th.~6.4]{pi01} and will be very useful
in the exact computation of the norm of certain relative Schur multipliers.

\begin{prp}\label{transfer}

  Let $0<p\le\infty$, let $\Gamma$ be a finite abelian group and
  $\Lambda\subseteq\Gamma$. Let $\row=\col=\Gamma$ and
  $I=\{(r,c)\in\Gamma\times\Gamma:r-c\in\Lambda\}$. 
  Let $G=\hat\Gamma$, so that $\Gamma$ is the group of characters on $G$. Let
  $\ell^p_{G,\Lambda}$ be the subspace of
  $\ell^p_G$ spanned by $\Lambda$. Then, for every
  sequence $\varphi$ on $\Lambda$, the norm of the relative Fourier
  multiplier $\varphi$ \(the norm of the linear operator
  $C_\varphi$ defined by $C_\varphi\gamma=\varphi_\gamma\gamma$\) on
  $\ell^p_{G,\Lambda}$ is bounded by the norm of the relative Schur
  multiplier $\breve\varphi=(\varphi_{r-c})_{(r,c)\in I}$ on $S^p_I$,
  and their c.b.~norms are equal.

\end{prp}

\begin{proof}
  For every $f\in \ell^p_{G,\Lambda}$, consider the multiplication
  operator $T_f$ on $\ell^2_G$ defined by $T_fh=fh$. $T_f$ is diagonal
  for the canonical basis $(\e_g)_{g\in G}$, so that
  $\|T_f\|_p=\|f\|_{\ell^p_{G}}$. In the basis $\Gamma$, its matrix
  representation is $T_f=\sum_{(r,c)\in I}\hat f(r-c)\e_{rc}$, so that
  $T_f\in S^p_I$ and $M_{\breve\varphi} T_f=T_{C_\varphi f}$. The norm
  of $\varphi$ on $\ell^p_{G,\Lambda}$ is therefore the norm of
  $\breve\varphi$ on the subspace of Toeplitz matrices in $S^p_I$.
  Furthermore, the c.b.~norm of $\varphi$ on $\ell^p_{G,\Lambda}$ is
  the norm of the Fourier multiplier $\varphi$ on
  $\ell^p_{G,\Lambda}(S^p)$ (\cite[Prop.~8.1.1]{pi98}).  For all
  $a_q\in S^p$ and $g$ in $G$, we have
\begin{equation*}
  \Bigl\|\sum_{q\in I}a_q\e_q\Bigr\|_p = \Bigl\|\sum_{(r,c)\in
    I}r(g)c(g)^{-1}a_{rc}\e_{rc}\Bigr\|_p
\end{equation*}
by matrix unconditionality (Eq.~\Ref{mu}), so that 
\begin{equation*}
\Bigl\|\sum_{q\in I}\breve\varphi_qa_q\e_q\Bigr\|_p 
  =\biggl\|\sum_{\gamma\in\Lambda}\varphi_\gamma\biggl(\sum_{\substack{(r,c)\in I\\
    r-c=\gamma}}a_{rc}\e_{rc}\biggr)\gamma\biggr\|_{\ell^p_{G}(S^p(S^p))}~.
\end{equation*}
It remains to note that $S^p(S^p)=S^p(\ell^2_\Gamma(\ell_2))$.
\end{proof}

\'Eric Ricard proposed the following device for computing the c.b.~norm of
a Schur multiplier, included here by his kind permission.

\begin{prp}\label{pos}

  Let $I=\N\times\N$ and $\varphi=(\varphi_q)_{q\in I}$ be the matrix
  representation of a positive operator. Then the c.b.~norm of the Schur
  multiplier $\varphi$ on $S^p$ is $\sup_n|\varphi_{nn}|$ for any $p\in\]0,\infty]$.
\end{prp}

\begin{proof}

  If $\varphi$ is positive, then $M_\varphi$ is a positive operator on
  $B(\ell^2)$. Henceforth its norm on $B(\ell^2)$ is attained at $\id$,
  and $\|\varphi\ast\id\|_\infty=\sup_n|\varphi_{nn}|$. This quantity is
  furthermore a lower bound on all $S^p$.
\end{proof}

\section{Idempotent Schur multipliers of norm one}\label{ss3}

A Schur multiplier is \textit{idempotent} if it is the indicator
function $\chi_I$ of some set $I\subseteq\row\times\col$; if $\chi_I$
is a Schur multiplier on $S^p$, then it is a projection of $S^p$ onto
$S^p_I$.  Idempotent Schur multipliers on $S^p$ and tensors in
$\ell^\infty_\col\projtens \ell^\infty_\row$ with $0,1$ coefficients
of norm $1$ may be characterized by the combinatorics of $I$.

\begin{prp}\label{complem}

  Let $I\subseteq\row\times\col$ be nonempty and $0<p\ne2<\infty$. The following
  are equivalent. 
  \begin{itemize}
  \item[$(a)$]For every finite rectangle set $\row'\times\col'$
    intersecting $I$
$$
    \biggl\| \sum_{(r,c)\in I\cap\row'\times\col'} \e_c\otimes\e_r
    \biggr\| _{\ell^\infty_\col\projtens \ell^\infty_\row}=1~.
$$
  \item[$(b)$]$S^p_I$ is completely $1$-complemented in $S^p$.
  \item[$(c)$]$S^p_I$ is $1$-complemented in $S^p$.
  \item[$(d)$]$I$ is a union of pairwise row and column disjoint
    rectangle sets, i.e.~complete bipartite graphs: there
    are pairwise disjoint sets $\row_j\subseteq\row$ and pairwise disjoint sets
    $\col_j\subseteq\col$ such that $I=\bigcup\row_j\times \col_j$.
  \end{itemize}

\end{prp}

\begin{proof} $(b)\imp(c)$ is trivial.\medskip

  $(a)\imp(b)$.  The c.b.~norm of a Schur multiplier $\varphi$ on
  $S^p$ is the supremum of the c.b.~norm of its restrictions
  $\varphi'=(\varphi_q)_{q\in\row'\times\col'}$ to finite rectangle
  sets $\row'\times \col'$. Furthermore, the c.b.~norm of an elementary
  Schur multiplier $(\eta_c\zeta_r)_{(r,c)\in\row\times\col}=\eta\otimes\zeta$ on
  $S^p$ equals
  $\|\eta\|_{\ell^\infty_\col}\|\zeta\|_{\ell^\infty_\row}$.\medskip

$(c)\imp(d)$. If $I$ is not a union of pairwise row and column disjoint rectangle sets,
then there are $r_0,r_1\in\row$, $c_0,c_1\in\col$ such that
$(r_0,c_0),\,(r_0,c_1),\,(r_1,c_0)\in I$ and $(r_1,c_1)\notin I$. Let $x(t)$,
$t\in\R$, be the operator from $\ell^2_{\{c_0,c_1\}}$ to $\ell^2_{\{r_0,r_1\}}$
with matrix $\Bigl(\begin{smallmatrix}1&\sqrt2\\ \sqrt2&t\end{smallmatrix}\Bigr)$. Its eigenvalues are
$$
\frac{1+t+\sqrt{9-2t+t^2}}2=2+\frac t3+o(t)~,\quad
\frac{1+t-\sqrt{9-2t+t^2}}2=-1+\frac{2t}3+o(t)~,
$$
so that
\begin{equation}\label{no4}
\|x(t)\|_\infty=2+\frac t3+o(t)~,~
\|x(t)\|_p^p=2^p+1+\frac p6(2^p-4)t+o(t)\quad{\rm for}~0<p<\infty
\end{equation}
and therefore $\|P_Ix(t)\|_p=\|x(0)\|_p>\|x(t)\|_p$ for some $t\ne0$ if
$p\ne2$.\medskip

$(d)\imp(a)$. Suppose $(d)$ and let $\row'\times\col'$ intersect $I$. Then
there are pairwise disjoint sets $\row_j$ and pairwise disjoint sets $\col_j$ such that $I\cap \row'\times
\col'=\row_1\times\col_1\cup\dots\cup\row_n\times\col_n$ and
$$
\sum_{(r,c)\in I\cap\row'\times\col'}\e_c\otimes\e_r = 
\sum_{j=1}^{n}\chi_{\col_j}\otimes\chi_{\row_j} = 
\mathop{\rm average}\limits_{\epsilon_j=\pm1}
\biggl(\sum_{j=1}^{n}\epsilon_j\chi_{\col_j}\biggr)
\otimes\biggl(\sum_{j=1}^{n}\epsilon_j\chi_{\row_j}\biggr)~
$$
which is an average of elementary tensors of norm $1$, so that its
projective tensor norm is bounded by $1$, and actually is equal to
$1$. 
\end{proof}

\begin{rem} Note that the proof of Prop.~\ref{complem} shows that the norm of a projection
$M_{\chi_I}\colon S^\infty\to S^\infty_I$ is either $1$ or at least $2/\sqrt3$,
as
$$
\left\|\begin{pmatrix}1&\sqrt2\\ \sqrt2&-1\end{pmatrix}
\right\|_\infty=\sqrt3~,\qquad
\left\|\begin{pmatrix}1&\sqrt2\\ \sqrt2&0\end{pmatrix}\right\|_\infty=2~.
$$
This is a non commutative analogue to the fact that an idempotent measure on a
locally compact abelian group $G$ has either norm $1$ or at least
$\sqrt5/2$ (\cite[Th.~3.7.2]{ru62}). $\|M_{\chi_I}\|$ actually equals
$2/\sqrt3$ for $I=\{(0,0),(0,1),(1,0)\}$, as shown in \cite[Lemma
3]{liv95}. In fact, the following decomposition holds:
\begin{multline*}
\e_0\otimes\e_0+\e_0\otimes\e_1+\e_1\otimes\e_0=\\
\bigl(
(\e^{-\iu\pi/12},\e^{\iu\pi/4})\otimes
(\e^{-\iu\pi/12},\e^{\iu\pi/4})+
(\e^{\iu\pi/12},\e^{-\iu\pi/4})\otimes
(\e^{\iu\pi/12},\e^{-\iu\pi/4})\bigr)/\sqrt3~.
\end{multline*}
\end{rem}

\begin{rem}
  The equivalence of $(c)$ with $(d)$ has been obtained
  independently by Banks and Harcharras \cite{bh04}.
\end{rem}

\section{Unconditional basic sequences}\label{ss4}

\begin{dfn}

Let $0<p\le\infty$ and $I\subseteq \row\times \col$. Let $\U=\T$ \(vs.~$\U=\{-1,1\}$\). 
\begin{itemize}
\item[$(a)$]$I$ is an \textit{unconditional basic sequence} in $S_p$ if there is a
  constant $D$ s.t.
\begin{equation}\label{u}
\Bigl\|\sum_{q\in I}\epsilon_qa_q\e_q\Bigr\|_p \leqslant D
\Bigl\|\sum_{q\in I} a_q\e_q\Bigr\|_p
\end{equation}
for every  choice of signs $\epsilon_q\in\U$
and every finitely supported sequence of complex coefficients $a_q$. Its
complex (\textit{vs.}~real) \textit{unconditional constant} is the
least such $D$.

\item[$(b)$]$I$ is a \textit{completely} unconditional basic sequence in $S_p$
  if there is a constant $D$ s.t.~\Ref{u} holds for every choice of signs
  $\epsilon_q\in\U$ and every finitely supported sequence of operator
  coefficients $a_q\in S^p$. Its complex (\textit{vs.}~real) \textit{complete}
    unconditional constant is the least such $D$.

\item[$(c)$]$I$ is a complex
\({\it vs.}~real, complex completely, real completely\)
\textit{$1$-un\-con\-di\-tio\-nal} basic sequence in $S_p$ if its complex
\({\it vs.}~real, complex complete, real complete\) unconditional
constant is $1$: Inequality \Ref{u} turns into the equality
\begin{equation}\label{1u}
\Bigl\|\sum_{q\in I}\epsilon_qa_q\e_q\Bigr\|_p =
\Bigl\|\sum_{q\in I} a_q\e_q\Bigr\|_p~.
\end{equation}

\end{itemize}

\end{dfn}

If Inequality~\Ref{u} holds for every real choice of signs, then it
also holds for every complex choice of signs at the cost of replacing
$D$ by $D\pi/2$ (see \cite{Se97}), so that there is no need to
distinguish between complex and real unconditional basic sequences.\medskip

The notions defined in $(a)$ and $(b)$ are called $\sigma(p)$ and
complete $\sigma(p)$ sets in \cite[\S4]{ha98} and \cite{hno01} (see
also the survey \cite[\S\,9]{px03}). The notions defined in $(c)$ are
their isometric counterparts.\medskip

By \cite[proof of Cor.~4]{sc78}, the real unconditional constant of
any basis of $S^p_I$ cannot be lower than a fourth of the real
unconditional constant of $I$ in $S^p$.\medskip

\begin{exa}
  A single column $\row\times\{c\}$, a single row $\{r\}\times\col$,
  the diagonal set ${\{({\rm row}\,n,\allowbreak{\rm
      col}\,n)\}}_{n\in\N}$ if $\row$ and $\col$ are copies of $\N$, are
  $1$-unconditional basic sequences in all $S^p$. In fact, every
  column section and every row section (this is the terminology of
  \cite[Def.~4.3]{va69}) is a $1$-unconditional basic sequence; note
  that the length of every path in the corresponding graph is at most 2.
\end{exa}

Note that the set $I$ is a (completely) $1$-unconditional basic
sequence in $S^p$ if and only if the set of relative Schur multipliers
by signs on $S^p_I$ define (complete) isometries. This yields by
Cor.~\ref{pinf}

\begin{prp}\label{pi}

  Let $I\subseteq\row\times\col$ and $0<p\leqslant\infty$. If $I$ is
  a real \({\it vs.}~complex\) $1$-unconditional basic sequence in
  $S^\infty$, then $I$ is also a real \({\it vs.}~complex\) completely
  $1$-unconditional basic sequence in $S^p$.

\end{prp}

\begin{exa}

  If $\row=\col=\{0,\dots,n-1\}$, $1\le p\le\infty$ and $I=\row\times
  \col$, then the complex unconditional constant of the basis of
  elementary matrices in $S^p$ is $n^{|1/2-1/p|}$ and coincides with
  its complete unconditional constant (see \cite[Lemma
  8.1.5]{pi98}). This is also the real unconditional constant if
  $n=2^k$ is a power of $2$ as it is then the norm of the Schur
  multiplier with matrix the $k$th tensor power
  $\bigl(\begin{smallmatrix}-1&1\\1&1\end{smallmatrix}\bigr) ^{\otimes
    k}$ (the $k$th Walsh matrix) on $S^p$. Let us now show that if $n=3$, the real
  unconditional constant of the basis of elementary matrices in
  $S^\infty$ is $5/3$ and differs from its complex unconditional
  constant, $\sqrt3$. In fact, because the canonical bases of
  $\ell^2_C$ and $\ell^2_R$ are symmetric, the norm of a Schur
  multiplier by real signs turns out to equal the norm of one of the three
  following Schur multipliers:
  $$
\begin{pmatrix}
  1&1&1\\1&1&1\\1&1&1
\end{pmatrix},
\begin{pmatrix}
  -1&1&1\\1&1&1\\1&1&1
\end{pmatrix},
\begin{pmatrix}
  -1&1&1\\1&-1&1\\1&1&-1
\end{pmatrix}.
$$
The first norm is $1$, the second one is $\sqrt2$ by comparison with
the Schur multiplier $\bigl(\begin{smallmatrix}-1&1\\1&1\end{smallmatrix}\bigr)$, and the
third one is by Prop.~\ref{transfer} the norm of the Fourier
multiplier $\varphi=(-1,1,1)$ on $\ell^\infty_G$ with $G$ the group of third roots
of unity: as this multiplier acts by convolution with $f=(1/3,-2/3,-2/3)$, its norm
is $\|f\|_{\ell^1_G}=5/3$. This may also be seen using
Prop.~\ref{pos}, as we have the following decomposition in positive
matrices:
$$
\begin{pmatrix}
  -1&1&1\\1&-1&1\\1&1&-1
\end{pmatrix} = 
\begin{pmatrix}
  1/3&1/3&1/3\\1/3&1/3&1/3\\1/3&1/3&1/3
\end{pmatrix}-
\begin{pmatrix}
  4/3&-2/3&-2/3\\-2/3&4/3&-2/3\\-2/3&-2/3&4/3
\end{pmatrix}.
$$
By complex interpolation, the real unconditional constant of the
  basis of elementary matrices is in fact strictly less than its
complex counterpart in all $S^p$.
\end{exa}

\section{\sloppy Varopoulos' characterization of unconditional matrices in $S^\infty$}\label{V}\label{ss5}

Our results may be seen as the isometric counterpart to results by
Varopoulos \cite{va69} on tensor algebras over discrete spaces and
their generalization to $S^p$. He characterized unconditional basic
sequences of elementary matrices in $S^\infty$ in his study of the
projective tensor product $c_0\projtens c_0$. We gather up his results
in the next theorem, as they are difficult to extract from the
literature.\medskip

\begin{thm}
  Let $I\subseteq\row\times\col$. The following are equivalent.
  \begin{itemize}
  \item[$(a)$]$I$ is an unconditional basic sequence in $S^\infty$.
  \item[$(b)$]$I$ is an interpolation set for Schur multipliers on $S^\infty$:
    every bounded sequence on $I$ is the restriction of a Schur multiplier on $S^\infty$.
  \item[$(c)$]$I$ is a $V$-Sidon set as defined in \cite[Def.~4.1]{va69}: every
    null sequence on $I$ is the restriction of the sequence of coefficients of
    a tensor in $c_0(\col)\projtens c_0(\row)$.
  \item[$(d)$]The coefficients of every tensor in
    $\ell^1_\col\injtens\ell^1_\row$ with support in
    $I$ form an absolutely convergent series.
  \item[$(e)$]$(z_cz_r)_ {(r,c)\in I}$ is a Sidon set in the dual of
    $\T^{\col\amalg\row}$, that is an unconditional basic sequence in
    $\mathscr{C}(\T^{\col\amalg\row})$.
  \item[$(f)$]There is a constant $\lambda$ such that for all
    $\row'\subseteq\row$ and $\col'\subseteq\col$ with $n$ elements
    $\#[I\cap\row'\times\col']\leqslant\lambda n$.
  \item[$(g)$]$I$ is a finite union of forests.
  \item[$(h)$]$I$ is a finite union of row sections and column sections.
  \item[$(i)$]Every bounded sequence supported by $I$ is a Schur
    multiplier on $S^\infty$.
\end{itemize}
\end{thm}
\begin{proof}[Sketch of proof]$(a)\imp(b)$. If $(a)$ holds, every sequence of signs
  $\epsilon\in\{-1,1\}^I$ is a Schur multiplier on $S^\infty_I$. By a
  convexity argument, this implies that every bounded sequence is a
  Schur multiplier on $S^\infty_I$, which may be extended to a Schur
  multiplier on $S^\infty$ with the same norm by
  \cite[Cor.~3.3]{pps89}.
\medskip 

$(b)\imp(c)$ holds by Grothendieck's inequality (see \cite[\S\,5]{pi01}) and an
approximation argument.\medskip 

$(d)$ is but the formulation dual to $(c)$ (see \cite[\S\,6.2]{va67}).\medskip

$(d)\imp(e)$. A computation yields
\begin{equation}\label{sidon}
\biggl\|\sum_{(r,c)\in I}\varphi_{rc}\e_c\otimes\e_r\biggr\|
_{\ell^1_\col\injtens\ell^1_\row} =
\sup_{z\in\T^{\col\amalg\row}}\biggl|\sum_{(r,c)\in I}\varphi_{rc}z_cz_r\biggr|~.
\end{equation}

$(e)\imp(f)$ is \cite[Th.~4.2]{va69}. (The proof can be found in
\cite[\S\,6.3]{va67} and in \cite[\S\,5]{va68}.)\medskip

$(f)\imp(g)$, $(f)\imp(h)$ can be found in \cite[Th.~6.1]{va68}.\medskip

$(g)\imp(h)$. In fact, a forest is a bisection in the terminology of
\cite[Def.~4.3]{va69}: it is the union of a row section $I_\row$
with a column section
$I_\col$. It suffices to
prove this for a tree: let the vertices of its edges be indexed by
words as described in the Terminology; then let $I_\row$ be the set of
its elements of the form $(v_{w\concat a},v_w)$ with $w$ of odd length
and $a$ a letter, and let $I_\col$ be the set of its elements of the
form $(v_w,v_{w\concat a})$ with $w$ of even length and $a$ a
letter.\medskip

$(h)\imp(i)$ is \cite[Th.~4.5]{va68}. Note that row sections and column sections
form $1$-unconditional basic sequences in $S^\infty$ and are $1$-complemented in
$S^\infty$ by Proposition \ref{complem}. \medskip

$(i)\imp(a)$ follows from the open mapping theorem.
\end{proof}

\section{Closed walk relations}\label{ss6}

We now introduce and study the combinatorial objects that we need in
order to analyze the expansion of the function defined by
\begin{equation}\label{Phidef}
\Phi_I(\epsilon,a)= 
\Bigl\|\sum_{q\in I}\epsilon_qa_q\e_q\Bigr\|_p^p
\end{equation}
for $I\subseteq\row\times\col$, a positive even integer $p=2k$, signs $\epsilon\in\T$ and
coefficients $a_q\in\C$, of which only a finite number are nonzero. In fact,
\begin{eqnarray*}
\Phi_I(\epsilon,a)&=&\trace\Bigl(\sum_{(r,c),(r',c')\in I}
(\epsilon_{rc}a_{rc}\e_{rc})^*
(\epsilon_{r'c'}a_{r'c'}\e_{r'c'})
\Bigr)^k\\
&=&\trace~\sum_{\substack{(r_1,c_1),(r'_1,c'_1),\dots,\\
(r_k,c_k),(r_k',c'_k)\in I}}
~\prod_{i=1}^k
(\epsilon_{r_ic_i}^{-1}\overline{a_{r_ic_i}}\e_{c_ir_i})
(\epsilon_{r'_ic'_i}a_{r'_ic'_i}\e_{r'_ic'_i})\\
&=&\sum_{\substack{(r_1,c_1),(r_1,c_2),\dots,\\
(r_k,c_k),(r_k,c_{k+1})\in I}}
~\prod_{i=1}^k\epsilon_{r_ic_i}^{-1}\epsilon_{r_ic_{i+1}}
\overline{a_{r_ic_i}}a_{r_ic_{i+1}}\quad(\text{where}~c_{k+1}=c_1)~.
\end{eqnarray*}
The latter sum runs over all closed walks $(c_1,r_1,c_2,\dots,c_k,r_k)$
of length $p$ in the graph $I$.  Its terms have the form
$\epsilon^{\beta-\alpha}\overline{a}^\alpha a^\beta$ in multinomial
notation. The attempt to describe those couples $(\alpha,\beta)$ that
effectively arise in this expansion yields the following definition.
\begin{dfn}\label{cwd}
Let $p=2k\geqslant0$ be an even integer and $I\subseteq\row\times\col$. 
\begin{itemize}
\item[$(i)$] Let $\Alpha_k^I=\{\alpha\in\N^I:\sum_{q\in
I}\alpha_q=k\}$ and set
$$
\textstyle
\P{k}^I=
\bigl\{(\alpha,\beta)\in\Alpha_k^I\times\Alpha_k^I:
\forall\, r\,\sum_c\alpha_{rc}=\sum_c\beta_{rc}~\text{and}~
\forall\, c\,\sum_r\alpha_{rc}=\sum_r\beta_{rc}\bigr\}.
$$
\item[$(ii)$] Two couples $(\alpha^1,\beta^1)\in\P{k_1}^I$,
  $(\alpha^2,\beta^2)\in\P{k_2}^I$ are \textit{row and column disjoint} if
  $k_1,k_2\geqslant1$ and 
$$\forall\,(r,c),(r',c)\in I\quad\alpha^1_{rc}\alpha^2_{r'c}=0~,\qquad
\forall\,(r,c),(r,c')\in I\quad\alpha^1_{rc}\alpha^2_{rc'}=0~.
$$
\item[$(iii)$] The set $\Pw{k}^I$ of \textit{closed walk relations} of length $p$ in $I$ is
  the subset of those $(\alpha,\beta)\in\P{k}^I$ that cannot be
  decomposed into the sum of two row and column disjoint couples.
\end{itemize}
\end{dfn}

\begin{exa} Let $\row=\col=\{0,1,2,3\}$ and $I=\row\times\col$. The couple
$(\e_{00}+\e_{11}+\allowbreak\e_{22}+\allowbreak\e_{33}, \e_{01}+\e_{10}+\e_{23}+\e_{32})$ is an
element of $\P{4}^I\setminus\Pw{4}^I$: it is the sum of the two row and
column disjoint closed walk relations
$(\e_{00}+\e_{11},\e_{01}+\e_{10}),(\e_{22}+\e_{33},\e_{23}+\e_{32})$.
\end{exa}

The next proposition shows that, for our purpose, 
closed walk relations describe entirely closed
walks.
\begin{prp}\label{gon:bijection}
Let $p=2k\geqslant0$ be an even integer and $I\subseteq\row\times\col$.  The mapping
$P=(c_1,r_1,c_2,r_2,\dots,c_k,r_k)\in(\col\times
\row)^k\mapsto(\alpha,\beta)\in\Alpha_k^I\times\Alpha_k^I$ with
\begin{eqnarray*}
\alpha_q&=&\card\bigl[i\in\{1,\dots,k\}:(r_i,c_i)=q\bigr]\\
 \beta_q&=&\card\bigl[i\in\{1,\dots,k\}:(r_i,c_{i+1})=q\bigr]\quad(\text{where}~c_{k+1}=c_1)
\end{eqnarray*}
is a surjection of the set of closed walks in the graph $I$ of length $p$
onto the set $\Pw{k}^I$ of  closed walk relations of
length $p$ in $I$. We shall write $P\sim (\alpha,\beta)$ and call
$n_{\alpha\beta}$ the number of closed walks of length $p$ mapped on
$(\alpha,\beta)$.
\end{prp}

\begin{proof} Let $(\alpha,\beta)\in\Pw{k}^I$. If $k=0$, the empty closed walk suits.
Suppose $k\geqslant1$; we have to find a closed walk of length $p$ that is mapped on
$(\alpha,\beta)$. Consider a walk
$(c_1,r_1,c_2,r_2,\dots,c_j,r_j,c_{j+1})$
in $I$ such that 
$\alpha^1_q=\card{\{i:(r_i,c_i)=q\}}\leqslant\alpha_q$ and 
$\beta^1_q=\card{\{i:(r_i,c_{i+1})=q\}}\leqslant\beta_q$ for every
$q\in\row\times\col$, and furthermore $j$ is maximal. 
We claim $(a)$ that $c_{j+1}=c_1$ and $(b)$ that $j=k$. Let
$(\alpha^2,\beta^2)=(\alpha,\beta)-(\alpha^1,\beta^1)$.\medskip

$(a)$. If $c_{j+1}\ne c_1$, then
$\sum_r\alpha^2_{rc_{j+1}}=\sum_r\beta^2_{rc_{j+1}}+1\geqslant1$. Thus there
is $r_{j+1}$ such that $\alpha^2_{r_{j+1}c_{j+1}}\geqslant1$. But then
$\sum_c\beta^2_{r_{j+1}c}=\sum_c\alpha^2_{r_{j+1}c}\geqslant1$ and there is
$c_{j+2}$ such that $\beta^2_{r_{j+1}c_{j+2}}\geqslant1$: $j$ is not
maximal.\medskip

$(b)$. Suppose $j<k$. Then $(\alpha^1,\beta^1)\in\P{j}^I$ and  
$(\alpha^2,\beta^2)\in\P{k-j}^I$. By hypothesis, they are not row and
column disjoint: there are $r,c,c'$ 
such that $\alpha^1_{rc}\alpha^2_{rc'}\geqslant1$ or $r,r',c$ such
that $\alpha^1_{rc}\alpha^2_{r'c}\geqslant1$. By interchanging $\row$ and
$\col$, by relabeling the vertices, we may suppose without loss of
generality that for $r'_1=r_j$ there is $c'_1$ such that
$\alpha^2_{r'_1c'_1}\geqslant1$. Then there is $c'_2$ such that
$\beta^2_{r'_1c'_2}\geqslant1$. By the argument used in Claim $(a)$, there is a
closed walk 
$(c'_1,r'_1,c'_2,\dots,c'_{j'},r'_{j'})$ such that
$\card{\{i:(r'_i,c'_i)=q\}}\leqslant\alpha^2_q$ and
$\card{\{i:(r'_i,c'_{i+1})=q\}}\leqslant\beta^2_q$ (with $c'_{j'+1}=c'_1$). 
Then the closed walk 
$$(c_1,r_1,c_2,r_2,\dots,c_j,r_j,c'_2,r'_2,
\dots,c'_{j'},r'_{j'},c'_1,r'_1)$$
shows that $j$ is not maximal.  
\end{proof}

\begin{exa} Let $I=\row\times\col=\{0,1\}\times\{0,1\}$. The closed walk relation
  $(\e_{00}+\e_{11},\allowbreak \e_{01}+\e_{10})\in\Pw{2}^I$ has two preimages by the surjection of
  Prop.~\ref{gon:bijection}: the two cycles $({\rm col}\,0,{\rm
    row}\,0,{\rm col}\,1,{\rm row}\,1)$ and $({\rm col}\,1,{\rm
    row}\,1,{\rm col}\,0,{\rm row}\,0)$.  The closed walk relation
  $(2\e_{00}+2\e_{01},\allowbreak 2\e_{00}+2\e_{01})$ has six
  preimages: closed walk $({\rm col}\,1,\allowbreak {\rm
    row}\,0,\allowbreak {\rm col}\,1,\allowbreak {\rm
    row}\,0,\allowbreak {\rm col}\,0,\allowbreak {\rm
    row}\,0,\allowbreak {\rm col}\,0,\allowbreak {\rm row}\,0)$ and
  three other like ones, $({\rm col}\,1,\allowbreak {\rm
    row}\,0,\allowbreak {\rm col}\,0,\allowbreak {\rm
    row}\,0,\allowbreak {\rm col}\,1,\allowbreak {\rm
    row}\,0,\allowbreak {\rm col}\,0,\allowbreak {\rm row}\,0)$ and
  another like one.
\end{exa}

We are now in position to state the following theorem, a matrix counterpart to the computation
presented in \cite[Prop.~2.5$(ii)$]{ne98}.
\begin{thm}\label{sigma:ldc}
Let $p=2k$ be a positive even integer and $I\subseteq\row\times\col$. 
\begin{itemize}
\item[$(a)$]The function $\Phi_I$ in Eq.~\Ref{Phidef} has the expansion
\begin{equation}\label{gon:expand}
\Phi_I(\epsilon,a) = 
\trace{\Bigl|\sum_{q\in I}\epsilon_qa_q\e_q\Bigr|^p} = \sum_{(\alpha,\beta)\in\Pw{k}^I}
n_{\alpha\beta}\epsilon^{\beta-\alpha}\overline{a}^\alpha
a^\beta~,
\end{equation}
where $n_{\alpha\beta}\ge1$ for every 
$(\alpha,\beta)\in\Pw{k}^I$.
\item[$(b)$] If the $a_q$ are coefficients in $S^p$, then
\begin{multline}\label{gon:expand2}
\Phi_I(\epsilon,a) = \trace\Bigl(\trace{\Bigl|\sum_{q\in
    I}\epsilon_qa_q\e_q\Bigr|^p}\Bigr)\\
= \sum_{(\alpha,\beta)\in\Pw{k}^I}
\epsilon^{\beta-\alpha}\sum_{(c_1,r_1,\dots,c_k,r_k)\sim(\alpha,\beta)}
\prod_{i=1}^ka_{r_ic_i}^*a_{r_ic_{i+1}}~({\rm with}~c_{k+1}=c_1)~.
\end{multline}
\end{itemize}
\end{thm}

\begin{proof} 
This follows from Def.~\ref{cwd} and Prop.~\ref{gon:bijection}.
\end{proof}

Note that the the edges of a closed walk $P\sim(\alpha,\beta)$
are precisely those $\{r,c\}$ such that $\alpha_{rc}+\beta_{rc}\ge1$.
$P$ is a cycle if and only if has not length $0$ or $2$ and
$\sum_r\alpha_{rc}\le1$ for all $c$ and
$\sum_c\alpha_{rc}\le1$ for all $r$. We now show how to decompose
closed walks into cycles.
\begin{prp}\label{gon:dec}
Let $P=(c_1,r_1,c_2,r_2,\dots,c_k,r_k)\sim(\alpha,\beta)$ be a closed walk.
\begin{itemize}
\item[$(a)$]If $r_i=r_j$ \({\it vs.}~$c_i=c_j$\) for some $i\ne j$,
  then $P$ is the juxtaposition of two nonempty closed walks
  $P_1\sim(\alpha^1,\beta^1)$ and $P_2\sim(\alpha^2,\beta^2)$ such
  that $(\alpha,\beta)=(\alpha^1,\beta^1)+(\alpha^2,\beta^2)$ and
  $\sum_c\alpha^1_{r_ic}\,,\,\sum_c\alpha^2_{r_ic}\ge1$ \({\it
    vs.}~$\sum_r\alpha^1_{rc_i}\,,\,\sum_r\alpha^2_{rc_i}\ge1$.\)
\item[$(b)$]$P$ is the juxtaposition of closed walks
  $P_j\sim(\alpha^j,\beta^j)$ such that
  $\sum_r\alpha^j_{rc}\in\{0,1\}$ for all $c$,
  $\sum_c\alpha^j_{rc}\in\{0,1\}$ for all $r$ and 
  $(\alpha,\beta)=\sum(\alpha^j,\beta^j)$.
\item[$(c)$]There are cycles $P_j\sim(\alpha^j,\beta^j)$ and $\gamma$ such that
$(\alpha,\beta)=(\gamma,\gamma)+\sum(\alpha^j,\beta^j)$.
\end{itemize}
\end{prp}

\begin{proof}
  $(a)$. If $r_i=r_j$ for $i<j$, we may suppose $j=k$: consider then
  $P_1=(c_1,\allowbreak r_1,\dots,\allowbreak c_i,\allowbreak r_i)$ and
  $P_2=(c_{i+1},r_{i+1},\dots,c_k,r_k)$. If $c_i=c_j$ for $i<j$, we
  may suppose $i=1$: consider then
  $P_1=(c_1,r_1,\dots,c_{j-1},r_{j-1})$ and
  $P_2=(c_j,r_j,\dots,c_k,r_k)$.\medskip
  
  $(b)$. Use $(a)$ in a maximality argument.\medskip
  
  $(c)$. Note that the closed walks $P_j$ in $(b)$ are either cycles
  or of length $2$; in the latter case $P_j=q\sim(\e_q,\e_q)$ for some
  $q\in I$.
\end{proof}

\section{Schur multipliers on a cycle}\label{ss7}
\label{sec:cst}

We can realize a cycle of even length $2s$, $s\ge2$, in the following
convenient way. Let $\Gamma=\Z/s\Z$. Then the adjacency relation of
integers modulo $s$ turns $\Gamma$ into the cycle $(0,1,\dots,s-1)$ of
length $s$. We double this cycle into the bipartite cycle $({\rm
  col}\,0,\allowbreak {\rm row}\,0,\allowbreak {\rm
  col}\,1,\allowbreak {\rm row}\,1,\dots,\allowbreak {\rm
  col}\,{s-1},\allowbreak {\rm row}\,{s-1})$ on $\Gamma\amalg
\Gamma$, corresponding to the set of couples $I=\{(i,i),(i,i+1):i\in
\Gamma\}\subseteq \Gamma\times \Gamma$.\medskip

$\Gamma$ is the group dual to $G=\hat \Gamma=\{z\in\C:z^s=1\}$.
We shall consider the space
$\ell^p_{G,\Lambda}$ spanned by $\Lambda=\{1,z\}$ in $\ell^p_{G}$,
where $z$ is the identical function on $G$: its
norm is given by
$\|a+bz\|_{\ell^p_{G}}=\big(\sum_{z^s=1}|a+bz|^p\bigr)^{1/p}$.\medskip

\begin{prp}\label{cyno}

  Let $\epsilon\in\T^I$ be a Schur multiplier by signs on $S^p_I$.
\begin{itemize}

\item[$(a)$]$\epsilon$ has the same norm as the Schur multiplier
  $\hat\epsilon$ given by $\hat\epsilon_q=1$ for
  $q\ne({s-1},0)$ and $\hat\epsilon_{s-1,0} =
  \overline{\epsilon_{00}} \epsilon_{01} \dots
  \overline{\epsilon_{{s-1},{s-1}}}
  \epsilon_{{s-1},0}$.

\item[$(b)$]$\epsilon$ has the same norm as $\check\epsilon$ given by
  $\check\epsilon_{ii} = 1$ and $\check\epsilon_{i,i+1}=\vartheta$
  with $\vartheta^s=\hat\epsilon_{{s-1},0}$.

\item[$(c)$]The norm of $\epsilon$ on $S^p_I$ is bounded below by the
  norm of the relative Fourier multiplier $\mu:a+bz\mapsto
  a+\vartheta bz$ on $\ell^p_{G,\Lambda}$; their c.b.~norms are equal.

\item[$(d)$]The norm of $\epsilon$ on $S^1_I$ and on $S^\infty_I$ is
  equal to the norm of $\mu$ on $\ell^1_{G,\Lambda}$ and on
  $\ell^\infty_{G,\Lambda}$: this norm is $\max_{z^s=-1}|\vartheta+z|/|1+\e^{\iu\pi/s}|$.

\item[$(e)$]$\epsilon$ defines an isometry on $S^p_I$ if and only if
  $\overline{\epsilon_{00}} \epsilon_{01} \dots
  \overline{\epsilon_{{s-1},{s-1}}}
  \epsilon_{{s-1},0}=1$ or $p/2\in\{1,2,\dots,s-1\}$.
\end{itemize}
\end{prp}

\begin{proof}

$(a)$ and $(b)$ follow from the 
matrix unconditionality of Schatten-von-Neu\-mann norms (see Eq.~\Ref{mu}). $(c)$
follows from Prop.~\ref{transfer}.
\medskip

$(d)$. Let $a=1$ if $s$ is even and $a=\e^{\iu\pi/s}$ if $s$ is
odd. Then $\|a+z\|_{\ell^1_{G}}=2\cot(\pi/2s)$ and $\|a+\vartheta
z\|_{\ell^1_{G}}=\max_{z^s=1}|\e^{\iu\pi/s}+\vartheta z|/\sin(\pi/2s)$ and this
yields a lower bound on the norm of $\mu$ on $\ell^1_{G,\Lambda}$.\medskip

The c.b.~norm of $\mu$ on $\ell^\infty_{G,\Lambda}$ is equal to its
norm and thus to the maximum of $f_\vartheta(r,t)=\|r\e^{\iu
  t}+\vartheta z\|_\infty / \|r\e^{\iu t}+z\|_\infty$ for $r\ge0$ and
$t\in[-\pi,\pi]$. As $f_\vartheta=f_{\omega\vartheta}$ if $\omega^s=1$
and $f_{\bar \vartheta}(r,t)=f_\vartheta(r,-t)$, we may suppose that
$u=\arg\vartheta\in[0,\pi/s]$. We may also restrict our study to
$t\in[u/2-\pi/s,u/2+\pi/s]$ and further to $t\in[u/2,u/2+\pi/s]$, as
$f_\vartheta(r,t)\le1$ if $t\in[u/2-\pi/s,u/2]$. In that case
$f_\vartheta(r,t)=|r\e^{\iu t}+\vartheta|/|r\e^{\iu\min(t,2\pi/s-t)}+1|$. It
turns out that $f_\vartheta(r,t)$ is maximal if $r=1$ and, as
$f_\vartheta(1,t)=\cos\bigl((u-t)/2\bigr)/ \cos\bigl(\min(t/2,\pi/s-t/2)\bigr)$
is increasing for $t\in[u/2,\pi/s]$ and decreasing for
$t\in[\pi/s,u/2+\pi/s]$, this maximum is $f_\vartheta(1,\pi/s)$.
\medskip

$(e)$. If $p$ is not an even integer and $\vartheta^s\ne1$, then
$\mu$ is not an isometry on $\ell^p_{G,\Lambda}$: otherwise the
functions $z$ and $\vartheta z$ would have the same distribution by
the Plotkin-Rudin Equimeasurability Theorem
\cite[Th.~I]{ru60}.\medskip

If $p\in\{2,4,\dots,2s-2\}$, then $I$ contains no cycle of length
$l\le p$, so that by Prop.~\ref{gon:dec}$(c)$, every closed walk
$P\sim(\alpha,\beta)$ satisfies $\alpha=\beta$ and the function
$\Phi_I(\epsilon,a)$ of Eq.~\Ref{gon:expand} is constant in
$\epsilon$.\medskip

If $p\in\{2s,2s+2,\dots\}$, the closed walk relation
$$(\alpha,\beta) = \Bigl(\sum_{i\in \Gamma}\e_{ii}, \sum_{i\in
  \Gamma}\e_{i,i+1}\Bigr)+(p/2-s)(\e_{00},\e_{00})$$
satisfies
$n_{\alpha\beta}\ge1$ by Prop.~\ref{gon:bijection}. Then the
coefficient of $\Phi_I(\epsilon,a)$ in ${\bar a}^\alpha a^\beta$
equals $n_{\alpha\beta}\overline{\epsilon_{00}} \epsilon_{01} \dots
\overline{\epsilon_{{s-1},{s-1}}}\epsilon_{{s-1},0}$ and must equal
the same quantity with $\epsilon$ replaced by $1$ if $\epsilon$
defines an isometry on $S^p_I$.
\end{proof}

\begin{rem} See \cite[p.~245]{li96} for a similar application of
  the Plotkin-Rudin Equimeasurability Theorem in $(e)$. 
\end{rem}

The real unconditional constant of $I$ is therefore the norm of
$\check\epsilon$ with $\vartheta=\e^{\iu\pi/s}$, and the complex
unconditional constant is the maximum of the norm of $\check\epsilon$ for
$\arg\vartheta\in[0,\pi/s]$ (note that the Schur multipliers $\epsilon$ and
$\bar\epsilon=(\overline{\epsilon_q})_{q\in I}$ have the same norm on $S^p_I$).
This yields
\begin{cor}\label{cyc}
  Let $0<p\leqslant\infty$ and $s\geqslant2$. Let $I$ be the cycle of
  length $2s$.
  \begin{itemize}
  \item[$(a)$]$I$ is a real $1$-unconditional basic sequence in $S^p$
    if and only if $p\in\{2,\allowbreak 4,\dots,\allowbreak 2s-2\}$. 
  \item[$(b)$]The real and complex unconditional constants of $I$
    in the spaces $S^1$ and $S^\infty$ equal $\sec\pi/2s$.
  \end{itemize}
\end{cor}

\section{1-unconditional matrices in $S^p$,  $p$ not an even integer}\label{ss8}

We now state the announced isometric counterpart to Varopoulos'
characterization of unconditional matrices in $S^\infty$ (Section \ref{V})
and its generalization to $S^p$ for $p$ not an even integer.
\begin{thm}\label{p-np}
  Let $I\subseteq\row\times\col$ be nonempty and
  $p\in\]0,\infty]\setminus2\N$.  The following are
  equivalent.
\begin{itemize}
\item[$(a)$]$I$ is a complex completely $1$-unconditional basic
  sequence in $S^p$.
\item[$(b)$]$I$ is a complex $1$-unconditional basic sequence in
  $S^p$.
\item[$(c)$]$I$ is a real $1$-unconditional basic sequence in $S^p$.
\item[$(d)$]$I$ is a forest.
\item[$(e)$]For each $\epsilon\in\T^I$ there are
 $\zeta\in\T^\col$  and $\eta\in\T^\row$ such that
  $\epsilon_{rc}=\zeta(c)\eta(r)$ for all $(r,c)\in I$.
\item[$(f)$]For each $\epsilon\in\{-1,1\}^I$ there are
 $\zeta\in\{-1,1\}^\col$  and $\eta\in\{-1,1\}^\row$ such that
  $\epsilon_{rc}=\zeta(c)\eta(r)$ for all $(r,c)\in I$.
\item[$(g)$]$I$ is a set of $V$-interpolation of constant $1$: for all
  $\varphi\in\ell^\infty_I$ 
\begin{equation}\label{interpol}
\inf\Bigl\{\biggl\|
\sum_{(r,c)\in \row\times\col}
\tilde{\varphi}_{rc}\e_c\otimes\e_r
\biggr\|_{\ell^\infty_\col\projtens \ell^\infty_\row}:
\tilde{\varphi}|_I=\varphi\Bigr\}=\sup_{q\in I}|\varphi_q|~.
\end{equation}
\item[$(h)$]$I$ is a $V$-Sidon set of constant $1$: for all
  $\varphi\in c_0(I)$ 
\begin{equation}\label{sid}
\inf\Bigl\{\biggl\|
\sum_{(r,c)\in \row\times\col}
\tilde{\varphi}_{rc}\e_c\otimes\e_r
\biggr\|_{c_0(\col)\projtens c_0(\row)}:
\tilde{\varphi}|_I=\varphi\Bigr\}=\sup_{q\in I}|\varphi_q|~.
\end{equation}
\item[$(i)$] For every tensor $u=\sum_{(r,c)\in
    I}\varphi_{rc}\e_c\otimes\e_r$ in
  $\ell^1_\col\injtens\ell^1_\row$ with support in $I$
  we have
  $\|u\|_{\ell^1_\col\injtens\ell^1_\row}=\sum_{(r,c)\in
    I}|\varphi_{rc}|$.
\item[$(j)$]$(z_cz_r)_ {(r,c)\in I}$ is a Sidon set of constant $1$ in
  the dual of $\T^{\col\amalg\row}$, that is a $1$-unconditional basic sequence
  in $\mathscr{C}(\T^{\col\amalg\row})$.
\item[$(k)$]For all $\row'\subseteq\row$ and $\col'\subseteq\col$ with
  $k\geqslant1$ elements $\#[I\cap\row'\times\col']\leqslant2k-1$.
\item[$(l)$]$I$ is an isometric interpolation set for Schur multipliers on
  $S^\infty$: every $\varphi\in\ell^\infty_I$ is the restriction of a
  Schur multiplier on $S^\infty$ with norm
  $\|M_\varphi\|=\|\varphi\|_{\ell^\infty_I}$.
\end{itemize}
\end{thm}

\begin{proof} $(a)\imp(b)\imp(c)$ is trivial.\medskip

$(c)\imp(d)$. Suppose that $I$ contains a cycle
$(c_0,r_0,\dots,c_{s-1},r_{s-1})$ with
$s\geqslant2$. Cor.~\ref{cyc}$(a)$ shows
that $I$ is not a real $1$-unconditional basic sequence in $S^p$.\medskip

$(d)\ssi(k)$. A tree on $2k$ vertices has exactly $2k-1$ edges, so
that a forest $I$ satisfies $(k)$. Conversely, a cycle of length $2s$
is a graph with $s$ row vertices, $s$ column vertices and $2s$
edges.\medskip

$(d)\imp(e)$. Let $I$ be a tree and index the vertices of its edges by words
$w\in W$ as described in the Terminology, so that $I$ is a set of form
$$\bigl\{(v_\emptyset,v_a),\,(v_{a\concat b},v_a),\,(v_{a\concat b},v_{a\concat
  b\concat c}),\,\dots\bigr\}\quad{\rm with}~\emptyset,a,a\concat b,a\concat
  b\concat c,\dots\in W~.$$
Define inductively $\eta$ and $\zeta$: let $\eta(v_\emptyset)=1$; suppose that
$\eta$ and $\zeta$ have been defined for all $v_w$ such that the length of $w\in
W$ is at most $2n$. We then set for all words $w$ of length $2n$ and all letters
$a,b$ such that $w,w\concat a,w\concat a\concat b\in W$:
$$\zeta(v_{w\concat a}) = \epsilon(v_w,v_{w\concat a}) / \eta(v_w)~,\quad
\eta(v_{w\concat a\concat b}) = \epsilon(v_{w\concat a\concat b}, v_{w\concat
  a}) / \zeta(v_{w\concat a})~.$$
If $I$ is a union of pairwise disjoint trees, we may define $\eta$ and $\zeta$ on each
tree separately; we may finally extend $\eta$ to $\row$ and $\zeta$ to $\col$ in an
arbitrary manner.\medskip
 
$(d)\imp(f)$ may be proved as $(d)\imp(e)$.\medskip

$(f)\imp(c)$. If $(f)$ holds, then every Schur multiplier by signs
$\epsilon\in\{-1,1\}^I$ is simple in the sense that $M_\epsilon x=D_\eta xD_\zeta$
with $D_\zeta$ ({\it vs.}~$D_\eta$) the diagonal operator on $\ell^2_\col$ ({\it
  vs.}~$\ell^2_\row$) of multiplication by $\zeta$ ({\it vs.}~$\eta$), so that the
c.b.~norm of $M_\epsilon$ on any $S^p_I$ is
$\|\zeta\|_{\ell^\infty_\col}\|\eta\|_{\ell^\infty_\row}=1$.\medskip

$(e)\imp(g)$. If $(e)$ holds, every
$\varphi\in\T^I\subseteq\ell^\infty_I$ may be extended to an
elementary tensor $\zeta\otimes\eta$ of norm $1$. $(g)$ follows
because every element of $\ell^\infty_I$ with norm $1$ is the half sum
of two elements of $\T^I$: note that $\e^{\iu t}\cos u =
\bigl(\e^{\iu(t+u)}+\e^{\iu(t-u)}\bigr)/2$.\medskip

$(g)\imp(h)$. It suffices to check Equality~\Ref{sid} for $\varphi$ with support
contained in a finite rectangle set $\row'\times \col'$. As $\ell^\infty_{\col'}\projtens \ell^\infty_{\row'}$
is $1$-complemented in $\ell^\infty_\col\projtens \ell^\infty_\row$,
Eq.~\Ref{interpol} yields Eq.~\Ref{sid}.\medskip

$(h)\ssi(i)$ because they are dual statements.\medskip

$(i)\ssi(j)$. Use Equality \Ref{sidon}.\medskip

$(h)\imp(l)$ may be deduced by the argument of
Prop.~\ref{complem}$(a)\imp(b)$.\medskip

$(l)\imp(a)$. Taking sign sequences $\varphi\in\T^I$ in $(l)$ shows that
all relative Schur multipliers by signs on
$S^\infty_I$ define isometries. Apply Prop.~\ref{pi}. \end{proof}
 
\begin{rem}
  The equivalence of $(e)$ with $(j)$ may also be shown as a consequence of the
  characterization of Sidon sets of constant $1$ in \cite{chm81}.
\end{rem}

Let us now answer Question~\ref{q3}.
\begin{cor}
  Let $I\subseteq\row\times\col$. The following are equivalent.
  \begin{itemize}
  \item[$(a)$]For all $\varphi\in c_0(I)$ one has $\bigl\|\sum_{(r,c)\in I}
    \varphi_{rc} \e_c \otimes \e_r
    \bigr\|_{c_0(\col)\projtens c_0(\row)}=\sup_{q\in I}|\varphi_q|$.
  \item[$(b)$]There are pairwise disjoint sets $\row_j\subseteq\row$ and pairwise disjoint sets
    $\col_j\subseteq\col$ such that $\row_j$ or $\col_j$ is a singleton for each $j$
    and $I=\bigcup \row_j\times \col_j$.
  \item[$(c)$]$I$ is a union of pairwise disjoint star graphs: every path in $I$ has
    length at most $2$.
  \end{itemize} 
\end{cor}

\begin{proof} $(a)\imp(b)$ follows from Prop.~\ref{complem}$(a)\imp(d)$ and
Th.~\ref{p-np}$(g)\imp(d)$.\medskip

$(b)\ssi(c)$. $(b)$ holds if and only if $(r,c),(r',c),(r,c')\in
I\imp(r=r'~{\rm or}~c=c')$ and therefore if and only if $(c)$ holds.\medskip

$(b)\imp(a)$. Suppose $(b)$ and let $\varphi\in c_0(I)$.
Let $\alpha_j=\sup_{(r,c)\in \row_j\times \col_j}|\varphi_{rc}|^{1/2}$. If
$\alpha_j=0$, let us define $\varrho^j=0$ and $\gamma^j=0$. Otherwise, if
$\row_j$ is a singleton $\{r\}$, let us define $\varrho^j=\alpha_j\e_r$ and
$\gamma^j$ by $\gamma^j_c=\varphi_{rc}/\alpha_j$ if $c\in \col_j$ and
$\gamma^j_c=0$ otherwise.  Otherwise, $\col_j$ is a singleton $\{c\}$ and we
define $\gamma^j=\alpha_j\e_c$ and $\varrho^j$ by
$\varrho^j_r=\varphi_{rc}/\alpha_j$ if $r\in \row_j$ and $\varrho^j_r=0$ otherwise.
Note that the $\gamma^j$ have pairwise disjoint support and are null sequences,
as well as the $\varrho^j$. Then
$$
\sum_{(r,c)\in I}\varphi_{rc}\e_c\otimes\e_r = 
\sum_j\gamma^j\otimes\varrho^j = 
\mathop{\rm average}\limits_{\epsilon_j=\pm1}
\biggl(\sum_j\epsilon_j\gamma^j\biggr)
\otimes\biggl(\sum_j\epsilon_j\varrho^j\biggr)~
$$
is an average of elementary tensors in $c_0(\col)\projtens c_0(\row)$
of norm $\sup_{q\in I}|\varphi_q|$, so that this average 
is also 
bounded by this norm, which obviously is a lower bound. \end{proof}

\section{1-unconditional matrices in $S^p$, $p$ an even integer}\label{ss9}

Let us now prove Theorem~\ref{th5} as a consequence of 
Theorem \ref{sigma:ldc} together with Proposition \ref{gon:dec}$(c)$.
\begin{thm}\label{sigma:THM}
Let $I\subseteq\row\times\col$ and $p=2k$ a positive even integer. The following
assertions are equivalent.
\begin{itemize}
\item[$(a)$]$I$ is a complex completely $1$-unconditional basic sequence in
$S^p$. 
\item[$(b)$]$I$ is a complex $1$-unconditional basic sequence in $S^p$.
\item[$(c)$]For every finite subset $F\subseteq I$ there is an operator $x\in
  S^p$, whose support $S$ contains $F$, such that
  $\bigl\|\sum\epsilon_qx_q\e_q\bigr\|_p$ does not depend on the complex choice
  of signs $\epsilon\in\T^{ S}$.
\item[$(d)$]$I$ is a real\index{real vs.~complex}\index{complex vs.~real}
$1$-unconditional basic sequence in $S^p$.
\item[$(e)$]For every finite subset $F\subseteq I$ there is an operator $x\in
  S^p$ with real matrix coefficients, whose support $S$ contains
  $F$, such that $\bigl\|\sum\epsilon_qx_q\e_q\bigr\|_p$ does not depend on the
  real choice of signs $\epsilon\in\{-1,1\}^{ S}$.
\item[$(f)$]Every closed walk $P\sim(\alpha,\beta)$ of length $2s\leqslant
  2k$ in $I$ satisfies $\alpha=\beta$.
\item[$(g)$]$I$ does not contain any cycle of length $2s\leqslant 2k$ as a subgraph.
\item[$(h)$]For each $v,w\in V$ there is at most one path in $I$ of
length $l\le k$ that joins $v$ to $w$.
\end{itemize}
\end{thm}

\begin{proof} $(a)\imp(b)\imp(c)$, $(b)\imp(d)\imp(e)$ are trivial.\medskip

$(c)\imp(g)$. Suppose that $I$ contains a cycle $P\sim(\gamma,\delta)$
of length $2s\leqslant 2k$: the corresponding set of couples is
$F=\{q:\gamma_q+\delta_q=1\}$. Let $x$ be as in $(c)$ and let
$(\alpha,\beta)=(\gamma,\delta)+(k-s)(\e_q,\e_q)$ for some arbitrary
$q\in F$.  Then $(\alpha,\beta)\in\Pw{k}^S$. Consider
$f(\epsilon)=\bigl\|\sum\epsilon_qx_q\e_q\bigr\|_p^p$ as a function on
the group $\T^S$. Then the Fourier coefficient of $f$ at the Steinhaus
character $\epsilon^{\beta-\alpha}$ is, by Th.~\ref{sigma:ldc}$(a)$,
\begin{eqnarray*}
&&\sum\bigl\{n_{\varepsilon\zeta}\overline{x}^\varepsilon x^\zeta:
(\varepsilon,\zeta)\in\Pw{k}^S\quad{\rm and}\quad
\zeta-\varepsilon=\beta-\alpha\bigr\}\\
&=&\overline{x}^\gamma x^\delta
\sum\bigl\{n_{\varepsilon\zeta}\overline{x}^{\varepsilon-\gamma}
x^{\zeta-\delta}:
(\varepsilon,\zeta)\in\Pw{k}^S\quad{\rm and}\quad
\zeta-\delta=\varepsilon-\gamma\bigr\}~.
\end{eqnarray*}
(Note that $\beta-\alpha=\delta-\gamma$.) As this last
sum has only positive terms and contains at least the term corresponding to
$(\alpha,\beta)$, $f$ cannot be constant.\medskip

$(e)\imp(g)$. Let $P\sim(\gamma,\delta)$, $F=\{q:\gamma_q+\delta_q=1\}$ and
$(\alpha,\beta)$ be as in the proof of the implication $(c)\imp(h)$. Let $x$ be
as in $(e)$.  Consider $f(\epsilon)=\bigl\|\sum\epsilon_qx_q\e_q\bigr\|_p^p$ as a
function on the group $\{-1,1\}^S$. Then the Fourier coefficient
$\widehat{f}(\epsilon^{\beta-\alpha})$ of $f$ at the
Walsh character $\epsilon^{\beta-\alpha}$ is, by Th.~\ref{sigma:ldc}$(a)$,
\begin{eqnarray*}
&&\sum\bigl\{n_{\varepsilon\zeta}x^{\varepsilon+\zeta}:
(\varepsilon,\zeta)\in\Pw{k}^S\quad{\rm and}\quad
\zeta-\varepsilon\equiv\beta-\alpha\pmod{2}\bigr\}\\
&=&x^{\gamma+\delta}
\sum\bigl\{n_{\varepsilon\zeta}x^{\varepsilon+\zeta-\gamma-\delta}:
(\varepsilon,\zeta)\in\Pw{k}^S\quad{\rm and}\quad
\zeta-\varepsilon\equiv\delta-\gamma\pmod{2}\bigr\}~.
\end{eqnarray*}
As this last sum has only positive terms and contains at least the term
corresponding to $(\alpha,\beta)$, $f$ cannot be constant.\medskip

$(f)\ssi(g)$. Apply Prop.~\ref{gon:dec}$(c)$.\medskip

$(g)\ssi(h)$. If $I$ contains a cycle $(v_0,\dots,v_{2s-1})$, then $I$
contains two distinct paths $(v_0,\dots,v_{s})$,
$(v_0,v_{2s-1},\dots,v_{s})$ of length $s$ from $v_0$ to $v_{s}$. If
$I$ contains two distinct paths $(v_0,\dots,v_l)$,
$(v'_0,\dots,v'_{l'})$ with $v_0=v'_0$, $v_{l}=v'_{l'}$ and $l,l'\le
k$, let $a$ be minimal such that $v_a\ne v'_a$, let $b\ge a$ be minimal such
that $v_b\in\{v'_a,\dots,v'_{l'}\}$ and let $d\ge a$ be minimal such that
$v'_d=v_b$. Then $(v_{a-1},\dots,v_b,v'_{d-1},\dots,v'_a)$ is a cycle
in $I$
of length $2s\le2k$.\medskip

$(f)\imp(a)$ holds by Theorem~\ref{sigma:ldc}$(b)$: If each
$(\alpha,\beta)\in\Pw{k}^I$ satisfies $\alpha=\beta$, then
Eq.~\Ref{gon:expand2} shows that 
$\Phi_I(\epsilon,z)$ as defined in
Eq.~\Ref{Phidef} is constant in $\epsilon$.
\end{proof}

\begin{rem}
  The equivalence $(b)\Leftrightarrow(g)$ is a non commutative
  analogue to \cite[Prop.~2.5$(ii)$]{ne98}.
\end{rem}

\begin{rem}
  In \cite[Th.~2.7]{ne99}, the condition of Th.~\ref{sigma:THM}$(f)$
  is visualized in another way: a closed walk
  $P=(c_1,r_1,\dots,c_s,r_s)\sim(\alpha,\beta)$ in $\N\times\N$ is
  considered as the polygonal closed curve $\gamma$ in $\C$ with sides
  parallel to the coordinate axes whose successive vertices are
  $r_1+\iu c_1$, $r_1+\iu c_2$, $r_2+\iu c_2$, \dots, $r_{s-1}+\iu
  c_s$, $r_s+\iu c_s$, $r_s+\iu c_1$ and again $r_1+\iu c_1$. Then
  $\alpha=\beta$ if and only if $\gamma$ is homologous (and homotopic)
  to zero with respect to the set of its points.
\end{rem}

\begin{rem} 
One cannot drop the assumption that $x$ has real matrix
coefficients in Th.~\ref{sigma:THM}$(e)$. Consider a $2\times2$
matrix $x$. Then $\trace x^*x=\sum|x_q|^2$ and $\det
x^*x=|x_{00}x_{11}-x_{01}x_{10}|^2$.  This shows that if
$\Re(\overline{x_{00}x_{11}}x_{01}x_{10})=0$, e.g.\
$x=\bigl(\begin{smallmatrix}1&1\\1&\iu\end{smallmatrix}\bigr)$, then the singular
values of $x$ do not depend on the real sign of the matrix
coefficients of $x$, whereas $({\rm col}\,0,{\rm row}\,0,{\rm
col}\,1,{\rm row}\,1)$ is a cycle of length $4$.
\end{rem}

\begin{rem}
  Theorem~\ref{sigma:THM}$(h)\imp(a)$ is the isometric counterpart to
  \cite[Th.~3.1]{hno01}. The following combinatorial problem arises
  naturally: if $I$ is such that the number of trails in $I$ of length $k$
  between two given vertices is uniformly bounded, is it
  so that $I$ is the union of a finite number of sets $I_j$ such that there
  is at most one path of length at most $k$ in $I_j$ between two given
  vertices~? In the simplest case, $k=2$, Ilijas Farah and Dominique
  Lecomte \cite{fl03} have deduced from \cite{ro85} that it is not so.
\end{rem}

\section{\sloppy Metric unconditional approximation pro\-per\-ty for $S^p_I$}\label{ss10}

Let $\row,\col$ be two copies of $\N$. It is well known that no $S^p$
has an unconditional basis or just a local unconditional structure
(see \cite[\S\,4]{px03}). $S^1$ and $S^\infty$ cannot even be
embedded in a space with unconditional basis. If $1<p<\infty$, then
$S^p$ has the unconditional finite dimensional decomposition
$$
\bigoplus_{n\in\N}S^p_{\{(r,c):r\le n,c=n\}}\oplus
S^p_{\{(r,c):r=n+1,c\le n\}}
$$
because the triangular projection associated to the idempotent Schur multiplier
$(\chi_{r\le c})$ is bounded on $S^p$.\medskip

\begin{dfn}
  Let $X$ be a separable Banach space and $\U=\T$ \(vs.~$\U=\{-1,1\}$\).
  \begin{itemize}
  \item [$(i)$]A sequence $(T_k)$ of operators on $X$ is an
    \textit{approximating sequence} if each $T_k$ has finite rank and
    $\|T_kx-x\|\to0$ for every $x\in X$. 
An approximating sequence of commuting projections
    is a \textit{finite-dimensional decomposition}.
  \item [$(ii)$](\cite{pw71}) The \textit{difference sequence}
    $(\Delta T_k)$ of $(T_k)$ is given by $\Delta T_1=T_1$ and $\Delta
    T_k=T_k-T_{k-1}$ for $k\ge2$. $X$ has the \textit{unconditional
      approximation property} $(uap)$ if there is an approximating
    sequence $(T_k)$ such that for some $D$
\begin{equation}\label{block:uapdef}
\biggl\|\sum_{k=1}^n\epsilon_k\Delta T_k\biggr\|
\leqslant D\quad\text{for all $n$ and $\epsilon_k\in\U$}~.
\end{equation}
The \textit{complex \(vs.~real\) unconditional constant} of $(T_k)$ is the least such $D$.
\item [$(iii)$](\cite[\S\,3]{ck91}, \cite[\S\,8]{gks93}) $X$ has the
  \textit{complex \(vs.~real\) metric unconditional approximation
    property} $(umap)$ if, for every $\delta>0$, it has an approximating
  sequence with complex \(vs.~real\) unconditional constant
  $1+\delta$. By \cite[Th.\ 3.8]{ck91} and \cite[Lemma 8.1]{gks93},
  this is the case if and only if there is an approximating sequence
  $(T_k)$ such that
\begin{equation}\label{block:umap:def}
\sup_{\epsilon\in\mathbb{S}}\|T_k+\epsilon (\id-T_k)\|\longrightarrow1~.
\end{equation}
\end{itemize}
\end{dfn}

$X$ has $(umap)$ if and
only if, for every given $\delta>0$, $X$ is isometric to a
$1$-complemented subspace of a space with a $(1+\delta)$-unconditional
finite-dimensional decomposition \cite[Cor.~IV.4]{gk97}. If $X$ has
$(umap)$, then, for any given $\delta>0$, $X$ is isometric to a
subspace of a space with a $(1+\delta)$-unconditional basis.
\medskip

\begin{exa}
  The simplest example is the subspace in $S^p$ of operators with an
  upper triangular matrix. In fact, if $I\subseteq\row\times\col$ is
  such that all columns $I\cap\row\times\{c\}$ ({\it vs.}~all rows
  $I\cap\{r\}\times\col$) are finite, then $S^p_I$ admits a
  $1$-unconditional finite-dimensional decomposition in the
  corresponding finitely supported idempotent Schur multipliers
  $\chi_{I\cap\row\times\{c\}}$ ({\it
    vs.}~$\chi_{I\cap\{r\}\times\col}$).
\end{exa}

Our results on complete $1$-unconditional basic sequences yield the
following theorem.
\begin{thm}
  Let $1\le p\le\infty$. Let $R_r\subset R$, $r\in\N$, be pairwise
  disjoint and finite. Let $C_c\subset C$, $c\in\N$, be pairwise
  disjoint and finite. Let $J\subset \N\times\N$ and
  $I=\bigcup_{(r,c)\in J}R_r\times C_c$. Then the sequence of Schur
  multipliers $(\chi_{R_r\times C_c})_{(r,c)\in J}$ forms a
  complex $1$-unconditional finite-dimensional
  decomposition for $S^p_I$ if and only if $J$ is a forest or $p$ is
  an even integer and $J$ contains no cycle of length $4,6,\dots,p$.
\end{thm}

We may always suppose that approximating sequences on spaces $S^p_I$
are associated to Schur multipliers. More precisely, we have
\begin{prp}
  Let $1\le p\le\infty$ and $I\subseteq \row\times \col$.  Let $(T_n)$ be an
  approximating sequence on $S^p_I$. Then there is a sequence of Schur
  multipliers $(\varphi_n)$ with finite support such that $(M_{\varphi_n})$ is
  an approximating sequence on $S^p_I$. If 
  $(T_n)$ satisfies \Ref{block:umap:def}, then so does
  $(M_{\varphi_n})$.
\end{prp}
\begin{proof}
  Let $\delta_n>0$ be such that $\delta_n\to0$.  As $T_n$ has finite
  rank, there is a finite $\row_n\times \col_n\subset \row\times \col$
  such that the projection $P_n$ of $S^p_I$ onto $S^p_{I\cap
    \row_n\times \col_n}$ satisfies $\|P_nT_n-T_n\|<\delta_n$. Let
  $U_n=P_nT_n$. Then $\|U_n+\epsilon(\id-U_n)\| <
  \|T_n+\epsilon(\id-T_n)\| + 2\delta_n$ for all
  $\epsilon\in\mathbb{S}$. For each $\zeta\in\mathbb{T}^{\col_n}$
  (\textit{vs.}~$\eta\in\mathbb{T}^{\row_n}$), let $D_\zeta$ ({\it
    vs.}~$D_\eta$) be the diagonal operator on $\ell^2_{\col_n}$ ({\it
    vs.}~$\ell^2_{\row_n}$) of multiplication by $\zeta$ ({\it
    vs.}~$\eta$): $D_\zeta$ and $D_\eta$ are unitary, so that
  $x\mapsto D_\eta xD_\zeta$ is an isometric operation on $S^p_{I\cap
    \row_n\times \col_n}$. Let $V_n$ be the operator defined by
  $$
  V_n(x)=\int_{\mathbb{T}^{\row_n}}{\rm d}\eta\int_{\mathbb{T}^{\col_n}}{\rm
    d}\zeta\, D_\eta^*U_n(D_\eta xD_\zeta)D_\zeta^*~.
  $$
  Then, as $U_ny\to y$ uniformly on compact sets,
  \begin{gather*}
  \|V_nx-x\|=\biggl\|\int_{\mathbb{T}^{\row_n}}{\rm
    d}\eta\int_{\mathbb{T}^{\col_n}}{\rm d}\zeta\, D_\eta^*\bigl(U_n(D_\eta
  xD_\zeta)-D_\eta xD_\zeta\bigr)D_\zeta^*\biggr\|\longrightarrow0~,\\
  \begin{split}
    \|V_nx+\epsilon(x-V_nx)\|&=\biggl\|\int_{\mathbb{T}^{\row_n}}{\rm
    d}\eta\int_{\mathbb{T}^{\col_n}}{\rm d}\zeta\, D_\eta^*\bigl((U_n+\epsilon(\id-U_n))(D_\eta xD_\zeta)\bigr)D_\zeta^*\biggr\|\\
  &\le\|U_n+\epsilon(\id-U_n)\|\,\|x\|\\
  &\le\bigl(\|T_n+\epsilon(\id-T_n)\|+2\delta_n\bigr)\|x\|~.
  \end{split}
  \end{gather*}
  A computation shows that $V_n$ is the operator associated to the Schur multiplier
  $\varphi_n={\bigl(\trace\e_q^*U_n(\e_q)\bigr)}_{q\in I}$ whose support is
  contained in ${I\cap \row_n\times \col_n}$.
\end{proof}

This proposition shows together with Cor.~\ref{pinf} the following results.
\begin{cor}
  Let $1\le p\le\infty$ and $I\subseteq \row\times \col$. 
  \begin{itemize}
  \item[$(a)$]If $S^p_I$ has $(umap)$, then  some
  sequence of Schur multipliers realizes it.
  \item[$(b)$]Let $J\subset I$. If $S^p_I$ has $(umap)$, then so does $S^p_J$.
  \item[$(c)$]If $S^\infty_I$ has $(umap)$, then so does $S^p_I$.
  \end{itemize}
\end{cor}

The arguments of \cite[\S\,6.2]{ne98} show \textit{mutatis mutandis}
\begin{thm}\label{umap}
  Let $1\le p\le\infty$, $I\subseteq \row\times \col$ and $\U=\T$
  \(vs.~$\U=\{-1,1\}$\). Consider the following properties.
  \begin{itemize}
  \item[$(a)$]$S^p_I$ has the property of $\tau$-unconditionality with
    $\tau$ the topology of pointwise convergence of matrix
    coefficients: for every $x\in S^p_I$ and every $\tau$-null
    sequence $(y_n)$
  $$
  \max_{\epsilon\in\mathbb{S}}\|x+\epsilon y_n\|_p -
  \min_{\epsilon\in\mathbb{S}}\|x+\epsilon y_n\|_p\longrightarrow0~.
  $$
  \item[$(b)$]$I$ enjoys the property $(\mathscr{U})$ of matrix block
  unconditionality in $S^p$: for each $\delta>0$ and finite $F\subset I$, there
  is a finite $G\subset I$ such that 
  $$\forall\,x\in B_{S^p_F}~\forall\,y\in B_{S^p_{I\setminus G}}\quad
    \max_{\epsilon\in\mathbb{S}}\|x+\epsilon y\|_p -
    \min_{\epsilon\in\mathbb{S}}\|x+\epsilon y\|_p < \delta~.
  $$
  \item[$(c)$]$S^p_I$ has $(umap)$.
  \end{itemize}
  Then $(c)\imp(a)\ssi(b)$. If $1<p<\infty$, then $(b)\ssi(c)$. If
  $p=1$, $S^1_I$ has $(umap)$ if and only if $S^1_I$ has $(uap)$ and
  $I$ enjoys $(\mathscr{U})$ in $S^1$.
\end{thm}

The case $p=\infty$ is extreme in the sense that the following properties are
equivalent for $S^\infty_I$: to be a dual space, to be reflexive, to
have a finite
cotype, not to contain $c_0$, because they are equivalent for I not to contain
any sequence $(r_n,c_n)$ with $(r_n)$ and $(c_n)$ injective, that
is for $I$ to be
contained in the union of a finite set of lines and a finite set of columns,
so that $S^\infty_I$ is isomorphic to $\ell^2_I$.\medskip

Let us now introduce the asymptotic property on $I$ that reflects the
combinatorics imposed by $(umap)$.
\begin{dfn}
  Let $I\subseteq\row\times\col$ and $k\ge1$. 
  \begin{itemize}
  \item[$(a)$]$I$ enjoys property $\mathscr{J}_k$ if for every path
    $P=(c_0,r_0,\dots,c_j,r_j)$ of odd length $2j+1\le k$ in $I$ there
    is a finite set $R'\times C'$ such that $P$ cannot be completed
    with edges in $I\setminus R'\times C'$ to a cycle of length
    $2s\in\{4j+2,\dots,2k\}$.
  \item[$(b)$]The \textit{asymptotic distance} $d_\infty(r,c)$ of
    $r\in\row$ and $c\in\col$ in $I$ is the supremum, over all finite
    rectangle sets $\row'\times\col'$, of the distance from $r$ to $c$
    in $I\setminus\row'\times\col'$.
  \end{itemize}

\end{dfn}

The asymptotic distance takes its values in $\{1,3,5,\dots,\infty\}$.
Note that $\mathscr{J}_k\imp\mathscr{J}_{k-1}$ and that
$\mathscr{J}_1$ is void. This implication is strict: let $\row,\col$
be two copies of $\N$ and, given $j\ge1$, consider the union $I_j$ of
all paths $(\text{col}\,0,\allowbreak \text{row}\,{nj+1},\allowbreak
\text{col}\,{nj+1},\allowbreak \dots,\allowbreak
\text{row}\,{nj+j},\allowbreak \text{col}\,{nj+j},\allowbreak
\text{row}\,0)$ of length $2j+1$. Then $I_j$ contains no cycle of
length $2s\in\{4,\dots,4j\}$ and therefore enjoys $\mathscr{J}_{2j}$,
but fails $\mathscr{J}_{2j+1}$;
$I_j\cup\{(\text{row}\,0,\text{col}\,0)\}$ contains no cycle of length
$2s\in\{4,\dots,2j\}$ and thus enjoys $\mathscr{J}_{j}$, but
fails $\mathscr{J}_{j+1}$. In particular, the properties
$\mathscr{J}_k$, $k\ge2$, are not stable under union with a
singleton. \medskip

Let us now explicit the relationship between $\mathscr{J}_k$ and $d_\infty$.

\begin{prp}
  Let $I\subseteq\row\times\col$ and $k\ge1$. 
  \begin{itemize}
  \item[$(a)$]$I$ enjoys $\mathscr{J}_k$ if and only if any two
    vertices $r\in R$ and $c\in C$ at distance $2j+1\le k$
    satisfy $d_\infty(r,c)\ge2k-2j+1$.
  \item[$(b)$]If $d_\infty(r,c)\ge2k+1$ for all $(r,c)\in\row\times\col$,
    then $I$ enjoys $\mathscr{J}_k$.
  \item[$(c)$]If $d_\infty(r,c)\le k$ for some $(r,c)\in\row\times\col$,
    then $I$ fails $\mathscr{J}_k$.
  \item[$(d)$]$I$ enjoys $\mathscr{J}_k$ for every $k$ if and only if
    $d_\infty(r,c)=\infty$ for every $(r,c)\in\row\times\col$.
  \end{itemize}
\end{prp}
\begin{proof}
  $(a)$ is but a reformulation of the definition of$\mathscr{J}_k$ and
  implies $(b)$. $(d)$ is a consequence of $(b)$ and $(c)$.\medskip

  $(c)$. If $d_\infty(r,c)\le k$, then there is $0\le j\le(k-1)/2$ such that
  there are infinitely many paths of length $2j+1$ from $c$ to $r$:
  there is a path
  $(c,r_1,c_1,\dots,r_j,c_j,r)$ that can be completed
  with edges outside any given finite set to a cycle of length
  $4j+2\le2k$.
\end{proof}
\begin{thm}
  Let $I\subseteq R\times C$ and $1\le p\le\infty$. If $p$ is an even
  integer then $S^p_I$ has complex or real $(umap)$ if and only if
  $I$ enjoys $\mathscr{J}_{p/2}$. If $p=\infty$ or if $p$ is not an even
  integer, then $S^p_I$ has real $(umap)$ only if
  $I$ enjoys $\mathscr{J}_k$ for every $k$.
\end{thm}

\begin{proof}
  Suppose that $I$ enjoys $(\mathscr{U})$ in $S^p$ and fails $\mathscr{J}_k$.
  Then, for some $s\le k$, $I$ contains a sequence of cycles $(c_0,\allowbreak
  r_0,\dots,\allowbreak c_{j-1},\allowbreak r_{j-1},\allowbreak
  c_j^n,\allowbreak r_j^n,\dots,\allowbreak c_{s-1}^n,\allowbreak
  r_{s-1}^n)$
  with the property that $\|x-y\|_p\le(1+1/n)\|x+y\|_p$ for all $x$ with support
  in $\{(r_0,\allowbreak c_0),\allowbreak (r_0,\allowbreak c_1),\dots,\allowbreak
  (r_{j-2},\allowbreak c_{j-1}),\allowbreak (r_{j-1},\allowbreak c_{j-1})\}$ and all $y$ with
  support in $\{(r_{j-1},\allowbreak c_j^n),\allowbreak
  (r_j^n,\allowbreak c_j^n),\dots,\allowbreak (r_{s-1}^n,\allowbreak c_{s-1}^n),\allowbreak
  (r_{s-1}^n,\allowbreak c_0)\}$. With the notation of Section~\ref{sec:cst}, this
  amounts to stating that the multiplier on
  $I=\{(i,i),(i,i+1)\}\subseteq\Z/s\Z\times\Z/s\Z$ given by
  $\epsilon_{rc}=1$ if $r,c\in\{0,\dots,j-1\}$ and $\epsilon_{rc}=-1$
  otherwise actually is an isometry on $S^p_I$. As
  $\overline{\epsilon_{00}} \epsilon_{01} \dots
  \overline{\epsilon_{{s-1}\,{s-1}}}
  \epsilon_{{s-1}\,0}=(-1)^{2s-2j+1}=-1$, this implies by
  Prop.~\ref{cyno}$(e)$ that $p/2\in\{1,2,\dots,s-1\}$.\medskip
  
  Suppose that $I$ enjoys $\mathscr{J}_k$. We claim that for every finite
  $F\subseteq I$ there is a finite $G\subseteq I$ such that every
  closed walk $P\sim(\alpha,\beta)$ of length $2k$ in $I$ satisfies
  $\sum_{q\in I\setminus G}\beta_q-\alpha_q=0$. That is, given a
  closed walk $(v_0,\dots,v_{2k-1})$ and
  $0=a_0<b_0<\dots<a_m<b_m<a_{m+1}=2k$ such that
  $v_{a_i},\dots,v_{b_i-1}\in I\setminus G$ and
  $v_{b_i},\dots,v_{a_{i+1}-1}\in F$,
  $$\bigl\{i\in\{0,\dots,m\}:a_i,b_i~\text{even}\bigr\} =
  \bigl\{i\in\{0,\dots,m\}:a_i,b_i~\text{odd}\bigr\}~.$$
  Suppose that this
  is not true: then there is an $s\le k$, there are
  $0=a_0<b_0<\dots<a_m<b_m<2s$ and there are cycles
  $(v_{a_0}^n,\dots,\allowbreak v_{b_0-1}^n,\allowbreak
  v_{b_0},\dots,\allowbreak v_{a_1-1},\dots,\allowbreak
  v_{a_m}^n,\dots,\allowbreak v_{b_m-1}^n,\allowbreak
  v_{b_m},\dots,\allowbreak v_{2s-1})$ such that the $(v_i^n)_{n\ge0}$
  are injective sequences of vertices and $b_i-a_i$ is even for at
  least one index $i$: let us suppose so for $i=0$. If $b_0-a_0\ge
  s-1$, consider the path
  $P=(v_{b_0},\dots,\allowbreak v_{a_0-1},\allowbreak
  v_{a_0}^0,\dots,\allowbreak v_{b_m-1}^0,\allowbreak
  v_{b_m},\dots,\allowbreak v_{2s-1})$
  of odd length $2s-1-(b_0-a_0)$; if $b_0-a_0\le s-1$, consider
  the path $P=(v_{2s-1},\allowbreak v_{a_0}^0,\dots,\allowbreak
  v_{b_0-1}^0,\allowbreak v_{b_0})$ of odd
  length $b_0-a_0+1$. Then $P$ can be completed with vertices
  outside any given finite set to a cycle of length at most $2s$
  because $(v_{2s-1},v_{a_0}^n,\dots,\allowbreak
  v_{b_0-1}^n,\allowbreak v_{b_0})$ is a path of length $b_0-a_0+1$ in
  $I$ for every $n$. This proves that $I$ fails $\mathscr{J}_s$.\medskip 
  
  The claim shows that $I$ enjoys $(\mathscr{U})$ in $S^p$ for $p=2k$. In
  fact, if $\tilde\epsilon\in\T^{F\cup(I\setminus G)}$ is defined by
  $\tilde\epsilon_q=1$ for $q\in F$ and
  $\tilde\epsilon_q=\epsilon\in\T$ for $q\in I\setminus G$, then, with
  the notation of Th.~\ref{sigma:ldc},
  $$\Phi_{F\cup(I\setminus G)}(\tilde\epsilon,a) =
  \sum_{(\alpha,\beta)\in\Pw{k}^{F\cup(I\setminus G)}}
  n_{\alpha\beta}\epsilon^{\sum_{q\in I\setminus G}
    \beta_q-\alpha_q}\overline{a}^\alpha a^\beta$$
  does not depend on
  $\epsilon$, so that $\|x+\epsilon y\|_{2k}=\|x+y\|_{2k}$ if $x\in
  S^{2k}_F$ and $y\in S^{2k}_{I\setminus G}$, and $S^{2k}_I$ has
  complex $(umap)$ by Th.~\ref{umap}$(b)\imp(c)$.
\end{proof}

\begin{rem}
  This theorem is a non commutative analogue to \cite[Th.~7.5]{ne98}.
\end{rem}

\section{Examples}\label{ss11}

One of Varopoulos' motivations for the study of the projective tensor
product $\ell_\infty\projtens\ell_\infty$ are lacunary sets in a
locally compact abelian group. Let $\Gamma$ be a discrete abelian
group and $\Lambda\subseteq\Gamma$. Let us say that $\Lambda$ is
\textit{$n$-independent} if every element of $\Gamma$ admits at most
one representation as the sum of $n$ elements of $\Lambda$, up to a
permutation. Let $G=\hat\Gamma$, so that $\Gamma$ is the group of
characters on $G$. Then the computation presented in
\cite[Prop.~2.5$(ii)$]{ne98} only in the case $\Gamma=\Z$ shows that
$\Lambda$ is a complex $1$-unconditional basic sequence in $L^{2n}(G)$
if and only if $\Lambda$ is $n$-independent. Furthermore $\Lambda$ is
a complex $1$-unconditional basic sequence in $L^{p}(G)$, $p$ not an
even integer, or in $\mathscr{C}(G)$, if and only if $\Lambda$ is
$n$-independent for all $n$. Real $1$-unconditional basic sequences
may be characterized in a similar way: the property replacing
$n$-independence is that every element of $\Gamma$ admits at most one
representation as the sum of $n$ elements of $\Lambda$, up to a
permutation and up to an even element. Consequently, if $\Gamma$
contains no element of order $2$, then real $1$-unconditional basic
sequences in $L^{p}(G)$ and in $\mathscr{C}(G)$ are also complex
$1$-unconditional.\medskip

Results on lacunary sets in $\Gamma$ transfer to lacunary
matrices in the following way.

\begin{prp}
  Let $\row,\col\subseteq\Gamma$ be such that each element of $\Gamma$ admits
  at most one representation as the sum of an element of $\row$ with an
  element of $\col$, so that $\Lambda\subseteq \row+\col$ determines uniquely
  the set $I\subseteq\row\times\col$ of all $(r,c)$ such that $r+c\in\Lambda$.
  \begin{itemize}
  \item[$(a)$]If $\Lambda$ is $n$-independent, then $I$ is a
    $1$-unconditional basic sequence in $S^{2n}$. 
  \item[$(b)$]The converse holds if furthermore each element of $\Gamma$ admits
    at most one representation as a sum of $n$ elements of $\row$ and $n$
    elements of $\col$, up to a permutation of the elements of $\row$ and a
    permutation of the elements of $\col$.
  \end{itemize}
\end{prp}
\begin{proof} $(a)$. If $\Lambda$ is $n$-independent, consider a closed walk
$P=(c_1,r_1,\dots,c_n,r_n)$ in $I$ with $P\sim\allowbreak(\alpha,\beta)$; as
$(r_1+c_1)+\dots+(r_n+c_n)=(r_1+c_2)+\dots+(r_{n-1}+c_n)+(r_n+c_1)$,
$(r_1+c_1,\ldots,r_n+c_n)$ is a permutation of
$(r_1+c_2,\ldots,r_{n-1}+c_n,r_n+c_1)$.
But then, by the hypothesis on $R$ and $C$,
$\bigl((r_1,c_1),\allowbreak\ldots,\allowbreak(r_n,c_n)\bigr)$ is a permutation
of
$\bigl((r_1,c_2),\allowbreak\ldots,\allowbreak(r_{n-1},c_n),\allowbreak(r_n,c_1)\bigr)$,
so that $\alpha=\beta$: $I$ is a complex completely $1$-unconditional
basic sequence by Th.~\ref{sigma:THM}$(f)\imp(a)$.\medskip

$(b)$. If $\Lambda$ is not $n$-independent, then there are two
$n$uples $(r_1+c_1,\allowbreak\ldots,\allowbreak r_n+c_n)$,
$(r'_1+c'_1,\allowbreak\ldots,\allowbreak r'_n+c'_n)\in I^n$ that are
not a permutation of each other such that
$(r_1+c_1)+\dots+(r_n+c_n)=(r'_1+c'_1)+\dots+ (r'_n+c'_n)$, so that if
we set $\alpha_q=\card{\{i:(r_i,c_i)=q\}}$ and
$\beta_q=\card{\{i:(r'_i,c'_i)=q\}}$ for each $q\in\row\times\col$,
then $\alpha\ne\beta$. But, by the hypothesis on $R$ and $C$,
$(r_1,\dots,r_n)$ is a permutation of $(r'_1,\dots,r'_n)$ and
$(c_1,\dots,c_n)$ is a permutation of $(c'_1,\dots,c'_n)$, so that
\begin{gather*}
\forall\,r~\sum\nolimits_c\alpha_{rc} = \card{\{i:r_i=r\}} =
\sum\nolimits_c\beta_{rc}~,\\
\forall\,c~\sum\nolimits_r\alpha_{rc} = \card{\{i:c_i=c\}} =
\sum\nolimits_r\beta_{rc}~,
\end{gather*}
and $(\alpha,\beta)\in\P{s}^I$. Then $(\alpha,\beta)$ is a sum of row
and column disjoint closed walk relations $(\alpha^j,\beta^j)$ such
that $\alpha^j\ne\beta^j$ for at least one $j$: $I$ is not a real
$1$-unconditional basic sequence in $S^{2n}$ by Th.~\ref{sigma:THM}$(d)\imp(f)$.\end{proof}

Harcharras \cite{ha98} used Peller's discovery \cite{pe82a} of
the link between Fourier and Hankel Schur multipliers to produce
unconditional basic sequences in $S^p$ that are unions of
antidiagonals in $\N\times\N$.  We have in our context the rather disappointing
\begin{prp}\label{hoory}
  Let $\Lambda\subseteq\N\subseteq\Z$ and
  $I=\{(r,c)\in\N\times\N:r+c\in\Lambda\}$.
  \begin{itemize}
  \item[$(a)$] $\Lambda$ is $2$-independent if and only if $I$ is a
    $1$-unconditional basic sequence in $S^4$.
    
  \item[$(b)$] If $\Lambda$ contains a set $\{i,j,k\}$ with three elements
    such that $i\leqslant j+k$, $j\leqslant i+k$ and $k\leqslant i+j$,
    then $I$ is a $1$-unconditional basic sequence in $S^p$ only if
    $p\in\{2,4\}$.
  \end{itemize}
\end{prp}

\begin{proof} $(a)$. $\Lambda$ is $2$-independent if and only if $i+l=j+k$ with
$i,j,k,l\in \Lambda$ has only trivial solutions. If this
is the case, consider a closed walk $P=(c_1,r_1,c_2,r_2)$ in $I$. As
$(r_1+c_1)+(r_2+c_2)=(r_1+c_2)+(r_2+c_1)$, we have $c_1=c_2$ or $r_1=r_2$, so
that $P$ is not a cycle. Conversely, if $\Lambda$ contains a solution
to $i+l=j+k$ with $i<j\leqslant k<l$, then $I$ contains the cycle $({\rm
  col}\,0,{\rm row}\,i,{\rm col}\,k-i,{\rm row}\,j)$.\medskip

$(b)$. We may suppose that $k>i,j$. Consider the cycle $({\rm col}\,0,{\rm
  row}\,i,{\rm col}\,k-i,\allowbreak {\rm row}\,i+j-k,{\rm col}\,k-j,{\rm row}\,j)$.
\end{proof}

Look up \cite[Def.~I.3.1]{bjl99} for the definition of a Steiner
system and \cite[Def.~1.3.1]{va98} for the definition of a 
generalized polygon.

\begin{prp}Let $2\le n\le m$, $I\subseteq R\times C$ with $\#C=n$ and $\#R=m$,
  and $e=\#I$.
  \begin{itemize}
  \item[$(a)$]If $I$ is a $1$-unconditional basic sequence in $S^4$,
    then 
    $$n\ge1 + \Bigl(\frac em-1\Bigr) + \Bigl(\frac
    em-1\Bigr)\Bigl(\frac en-1\Bigr)~,$$
    that is
    $e^2-me-mn(n-1)\le0$. Equality holds if and only if $I$ is the
    incidence graph of a Steiner system $S(2,e/m;n)$ on $n$ points and
    $m$ blocks.
  \item[$(b)$]If $I$ is a $1$-unconditional basic sequence in $S^6$,
    then 
    $$n\ge1 + \Bigl(\frac em-1\Bigr) + \Bigl(\frac
    em-1\Bigr)\Bigl(\frac en-1\Bigr) + \Bigl(\frac
    em-1\Bigr)^2\Bigl(\frac en-1\Bigr)~,$$
    that is
    $e^3-(m+n)e^2+2mne-m^2n^2\le0$. Equality holds if and only if $I$ is
    the incidence graph of the quadrangle \(the cycle of length $8$\) or
    of a generalized quadrangle with $n$ points and $m$ lines.
  \item[$(c)$]If $I$ is a $1$-unconditional basic sequence in $S^{2k}$
    with $k\ge1$ an integer, then 
    \begin{equation}\label{hoo}
n\ge\sum_{i=0}^k\Bigl(\frac em-1\Bigr)^{\lceil\frac i2\rceil}\Bigl(\frac
 en-1\Bigr)^{\lfloor\frac i2\rfloor}~.
\end{equation}
Equality holds if $I$ is the incidence graph of the $(k+1)$-gon \(the
cycle of length $2k+2$\) or of a generalized $(k+1)$-gon with $n$
points and $m$ lines.
\end{itemize}
\end{prp}

\begin{proof}
  By Theorem~\ref{sigma:THM}$(b)\imp(g)$, $I$ is a $1$-unconditional basic
  sequence in $S^{2k}$, with $k\ge1$ an integer, if and only if $I$ is a
  graph of girth $2k+2$ in the sense of \cite{ho02}. Therefore $(a)$ and $(b)$ are shown in
  \cite[Prop.~4, Th.~8, Rem.~10]{ne01}. Inequality \Ref{hoo} is
  \cite[Eq.~(1)]{ho02} and the sufficient condition for equality follows from 
  \cite[Lemma 1.5.4]{va98}. 
\end{proof}

Consult \cite[Tables~A1.1, A5.1]{bjl99} for examples of Steiner
systems and \cite[Table~2.1]{va98} for examples of generalized
polygons. In both cases, the corresponding incidence graph is biregular: every
vertex in $R$ has equal degree $s+1$ and every vertex in $C$ has equal
degree $t+1$. Arbitrarily large generalized $(k+1)$-gons exist only
if $2k\in\{4,6,10,14\}$ (\cite[Lemma 1.7.1]{va98}); for
$2k\in\{6,10,14\}$, it follows from \cite[Lemma 1.5.4]{va98} that 
$$n=(s+1)\frac{(st)^{(k+1)/2}-1}{st-1}~,~
   m=(t+1)\frac{(st)^{(k+1)/2}-1}{st-1}~.$$

\begin{rem}
  Let $I\subseteq R\times C$ with $\#C=\#R=n$. Inequality~\Ref{hoo} shows
  that if $I$ is a $1$-unconditional basic sequence in $S^{2k}$, then
  $\#I\le n^{1+1/k}+(s-1)n/s$. If $p\notin\{4,6,10\}$, the existence
  of $1$-unconditional basic sequences in $S^{2k}$ such that
  $\#I\succcurlyeq n^{1+1/k}$ is in fact an important open problem in Graph
  Theory.
\end{rem}

Laboratoire de Math\'ematiques\\
Universit\'e de Franche-Comt\'e\\
25030 Besan\c con cedex, France\\
neuwirth@math.univ-fcomte.fr\medskip

Key words: Schatten-von-Neumann class, $1$-unconditionality, $V$-Sidon
set, cycle in a graph, metric unconditional approximation property. \medskip

2000 Mathematics Subject Classification: 47B10, 46B15, 46B04, 43A46,
05C38, 46B28.

\end{document}